\magnification=1200 
\parskip = 10pt plus 5pt 
\parindent = 15pt 
\baselineskip = 15pt 
\input amssym.def 
\input amssym
\input epsf.tex
\def\ssquare{\vrule width.6em height.5em depth.1em\relax}
\def\qed{\ifhmode\unskip\nobreak\fi\quad\ssquare}
\line{\hfil{November, 2000~~~~~~}}
\topglue 1.05in 
\pageno = 0
\footline={\ifnum \pageno < 1 \else \hss \folio \hss \fi} 
\centerline{\bf A RATIONALITY CONDITION FOR THE EXISTENCE}
\centerline{\bf OF ODD PERFECT NUMBERS}
\vskip .15in 
\centerline{Simon Davis}
\vskip .15in 
\centerline{School of Mathematics and Statistics}
\vskip 1pt
\centerline{University of Sydney}
\vskip 1pt
\centerline{NSW 2006, Australia}
\vskip 20pt 
\noindent{\bf Abstract.}  A rationality condition is derived for the 
existence of odd perfect numbers involving the square root of a product, 
which consists of a sequence of repunits multiplied by twice the
base of one of the repunits.   This constraint also provides 
an upper bound for the density of odd integers 
which could satisfy ${{\sigma(N)}\over N}=2$, where $N$ belongs to a fixed
interval with a lower limit greater than $10^{300}$.
Characteristics of prime divisors of repunits are used to 
establish whether the product containing the repunits can 
be a perfect square.   It is shown that the arithmetic 
primitive factors of the repunits with different prime bases can be equal 
only when the exponents are different, with possible exceptions 
derived from solutions of a prime equation.  This equation is one example of 
a more general prime equation, ${{q_j^n-1}\over {q_i^n-1}}=p^h$ and the 
demonstration of the non-existence of solutions when $h\ge 2$ requires the 
proof of a special case of Catalan's conjecture.   Results concerning the 
exponents of prime divisors of the repunits are obtained, and they are 
combined with the method of induction to prove a general theorem on the 
non-existence of prime divisors satisfying the rationality condition.
\vskip .30in
\noindent

\vskip .15in


\vfill\eject
\noindent
{\bf 1. Introduction}

The algorithm for demonstrating the non-existence of odd perfect numbers
with fewer than nine different prime divisors requires the expansion
of the ratio ${{\sigma(N)}\over N}$ and strict inequalities imposed on the
sums of powers of the reciprocal of each prime divisor [1][2].
Although it is possible to establish that ${{\sigma(N)}\over N}\ne 2$ 
when $N$ is divisible by certain primes, there are odd integers with a given 
number of prime divisors such that ${{\sigma(N)}\over N}> 2$ , 
while ${{\sigma(N)}\over N} < 2$ for other integers with the same number of
distinct prime factors. Moreover, 
the range of the inequality for $\left\vert {{\sigma(N)}\over N} - 
2\right\vert$ can be made very small even when $N$ has a few prime factors.  
Examples of odd integers with only five 
distinct prime factors have been found that produce a ratio nearly equal 
to 2:  $\left\vert {{\sigma(N)}\over N} - 2\right\vert < 10^{-12}$ [3]. 
Since it becomes progressively more difficult to establish the inequalities
as the number of prime factors increases, a proof by method of induction
based on this algorithm cannot be easily constructed.

In $\S 2$, it will be shown that there is a rationality condition for
the existence of odd perfect numbers.  Setting ${{\sigma(N)}\over N}$ 
equal to 2 is equivalent to equating the square root of a product, 
$2(4k+1)\prod_{i=1}^\ell~{{q_i^{n_i}-1}\over {q_i-1}}~
{{(4k+1)^{4m+2}-1}\over {4k}}$, which contains a sequence of repunits, with
a rational number.   This relation provides both an upper bound for the 
density of odd perfect numbers in any fixed interval in ${\Bbb N}$ with 
a lower limit greater than $10^{300}$ and a direct analytical method 
for verifying their non-existence, since it is based on the irrationality 
of the square root of any unmatched prime divisors in the product.  
This condition is used in $\S 3$ to demonstrate the non-existence of a 
special category of odd perfect numbers.  The properties of prime divisors 
of Lucas sequences required for the study of the square root of the product 
of the repunits are described in $\S 4$ and $\S 5$.  An induction argument 
is constructed in $\S 6$, which proves that the square root expression is not 
rational for generic sets of prime divisors, each containing a large 
number of elements.  This is first established for odd integers with four 
distinct prime divisors and then by induction using the properties of the 
divisors of the repunits.

\vskip 10pt
\noindent
{\bf 2. Rationality Condition for the Existence of Odd Perfect Numbers}

From the condition for the odd integer $N=(4k+1)^{4m+1}s^2=(4k+1)^{4m+1}
q_1^{2\alpha_1}...q_\ell^{2\alpha_\ell}$, $gcd(4k+1,s)=1$ [4]-[6] to be a 
perfect number, 
$$\eqalign{{\sigma(N)}\over N}~=~\left[{{(4k+1)^{4m+2}-1}
\over {4k(4k+1)^{4m+1}}}\right]{{\sigma(s^2)}\over {s^2}}
~=~\left[{{(4k+1)^{4m+2}-1}\over
{4k(4k+1)^{4m+1}}}\right]\left[{{\sigma(s)^2}\over {s^2}}\right]
\left[{{\sigma(s^2)}\over {\sigma(s)^2}}\right]~=~2
\eqno(1)
$$
it follows that
$${{\sigma(s)}\over s}~=~{\sqrt 2}~\prod_{i=1}^\ell~{{(q_i^{\alpha+1}-1)}
\over {(q_i~-~1)^{1\over 2}(q_i^{2\alpha_i+1}~-~1)^{1\over 2}}}
~\times~\left[{{4k(4k+1)^{4m+1}}\over {(4k+1)^{4m+2}~-~1}}\right]^{1\over 2}
\eqno(2)
$$ 
and
$$\prod_{i=1}^\ell~{1\over {(q_i^{\alpha_i+1}~-~1)}}~{{\sigma(s)}\over s}
~=~{\sqrt 2}~\prod_{i=1}^{\ell}~{1\over {(q_i^{2\alpha_i+1}~-~1)^{1\over 2}
(q_i~-~1)^{1\over 2} }}~\times~ \left[{{4k(4k+1)^{4m+1}}\over {(4k+1)^{4m+2}
~-~1}}\right]^{1\over 2}
\eqno(3)
$$
Consistency of equation (3) for finite $\ell$ requires rationality of the 
entire square root expression.

The known integer solutions to ${{x^n~-~1}\over {x~-~1}}~=~y^2$ [7]-[9]
do not include the pairs $(x,n)=(4k+1,4m+2)$, implying that
$\left[{{(4k+1)^{4m+2}~-~1}\over {4k}}\right]^{1\over 2}$ is not a rational 
number.  The number $[1+q_i+q_i^2+...+q_i^{2\alpha_i}]^{1\over 2}$ is 
only rational when $q_i=3,~\alpha_i=2$, so that if $3$ is a prime factor 
of $s$

$$\eqalign{\prod_{i=1}^\ell~{1\over {(q_i^{2\alpha_i+1}~-~1)^{1\over 2}}}
{1\over {(q_i-1)^{1\over 2}}}~&=~\prod_{i=1}^\ell~{1\over 
{(q_i^{2\alpha_i+1}-1)}}~[ 1~+~q_i~+~q_i^2~+~...~+~q_i^{2\alpha_i}]^{1\over 2}
\cr
~=~\left({{11}\over {242}}\right)^{\delta_{q_i,3}\delta_{\alpha_i,2}}
~\times~&\prod_{{i=1}\atop {(q_i,2\alpha_i+1)\ne (3,5)}}^\ell
~{1\over {(q_i^{2\alpha_i+1}-1)} }~
[1~+~q_i~+~q_i^2~+~...~+~q_i^{2\alpha_i}]^{1\over 2}
\cr}
\eqno(4)
$$

From equation (3), the non-existence of odd perfect numbers can be deduced if

$$\eqalign{2(4k+1)
\prod_i~&(1~+~q_i~+~q_i^2~+~...~q_i^{2\alpha_i})
\cr
&~~~~\cdot~(1~+~(4k+1)~+~(4k+1)^2~+~...~(4k+1)^{(4m+1)})
\cr
~&~~~~~~~~=~2(4k+1)\prod_i~{{q_i^{2\alpha_i+1}~-~1}\over {q_i~-~1}}
~{{(4k+1)^{4m+2}~-~1}\over {4k}}
\cr}
\eqno(5)
$$
is not the square of an integer, with $q_i\ne 3$ or $\alpha_i\ne 2$.

Since the repunit ${{x^n~-~1}\over {x~-~1}}$ is the Lucas sequence
$$U_n(a,b)~=~{{\alpha^n~-~\beta^n}\over {\alpha~-~\beta}}
\eqno(6)
$$
with $\alpha=x$, $\beta=1$, derived from the second-order recurrence 
relation 
$$U_{n+2}(a,b)~=~a~U_{n+1}(a,b)~-~b~U_n(a,b)
\eqno(7)
$$
where $a=\alpha+\beta=x+1$ and $b=\alpha\beta=x$, 
the rationality condition can be applied equally well to the product
$\left[2(4k+1)\prod_{i=1}^\ell~U_{2\alpha_i+1}(q_i+1,q_i)\cdot 
U_{4m+2}(4k+2,4k+1)\right]^{1\over 2}$.

The number of square-full integers up to
$N$ is $N^{1\over 2}-{3\over 2}N^{-1}+O(N^{-{3\over 2}})$.  With a lower bound
of $10^{300}$ for an odd perfect number [10], it follows that
$2(4k+1)\cdot\prod_{i=1}^\ell~(q_i^{2\alpha_i}+O(q_i^{2\alpha_i-1}))
\cdot ((4k+1)^{4m+1}+O((4k+1)^{4m})) > 10^{301}$.  Given a lower bound
of $10^6$ for the largest prime factor [11], $10^4$ for the second largest 
prime factor and $10^2$ for the third largest prime factor of $N$ [12], 
the density of prime products $(4k+1)\prod_{i=1}^\ell~q_i$, 
given by $\prod_{i=1}^\ell~{1\over {ln~q_i}}\times {1\over {ln~(4k+1)}}$, 
is bounded above by $8.032\times 10^{-5}$ when there are eight different prime 
factors [2] and $1.004\times 10^{-6}$ when there are eleven 
different prime factors not including 3 [13][14].  Given that the
probability of an integer being a square is independent of it being expressible
in terms of a product of repunits, the density of square-full numbers having
the form $2(4k+1)~\sigma\left((4k+1)^{4m+1}\prod_i~q_i^{2\alpha_i}\right)$ 
in the interval $[N^\ast, N^\ast + N_0]$, where $N^\ast > 10^{301}$ and 
$N_0$ is a fixed number, is bounded above by $3.28\times 10^{-159}$ when 
there are at least eight different prime factors and $5.13\times 10^{-163}$ 
when $N$ is relatively prime to 3 and has more than ten different prime 
factors.  
 
Although neither the numerator or the denominator will be 
squares of integers when $q_i\ne 3$ or $\alpha_i\ne 2$, there still 
remains the possibility that the terms could be equal multiples of 
different squares.  Since the repunit ${{x^n~-~1}\over {x~-~1}}$ 
is the Lucas sequence
$$U_n(a,b)~=~{{\alpha^n~-~\beta^n}\over {\alpha~-~\beta}}
\eqno(8)
$$
with $\alpha=x$, $\beta=1$, derived from the second-order recurrence  
$$U_{n+2}(a,b)~=~a~U_{n+1}(a,b)~-~b~U_n(a,b)
\eqno(9)
$$
where $a=\alpha+\beta=x+1$ and $b=\alpha\beta=x$, 
the rationality condition can be applied to the
$\left[2(4k+1)\prod_{i=1}^\ell~U_{2\alpha_i+1}(q_i+1,q_i)\cdot 
U_{4m+2}(4k+2,4k+1)\right]^{1\over 2}$.

\vfill\eject
\noindent
{\bf 3. Proof of the non-existence of odd perfect numbers for a special class}
\hfil\break 
\phantom{.....}{\bf of integers}

The even repunit ${{(4k+1)^{4m+2}-1}\over {4k}}$ contains only a single
power of 2 since $1+(4k+1)+(4k+1)^2+...+(4k+1)^{4m+1}\equiv 4m+2
\equiv 2~(mod~4)$.  Thus, the rationality condition can be applied to
a product of odd numbers $\left[(4k+1)\prod_{i=1}^\ell
U_{2\alpha_i+1}(q_i+1,q_i){1\over 2}U_{4m+2}(4k+2,4k+1)\right]^{1/2}$.  
Suppose
$$\eqalign{\prod_{i=1}^\ell~{{q_i^{2\alpha_i+1}-1}\over {q_i-1}}\cdot 
\left[{{8k(4k+1)}\over {(4k+1)^{4m+2}-1}}\right]~&=~{{r^2}\over {t^2}}
\cr
\prod_{i=1}^\ell~{{q_i^{2\alpha_i+1}-1}\over {q_i-1}}~(4k+1)~t^2
~&=~{{(4k+1)^{4m+2}-1}\over {8k}}~r^2
\cr}
\eqno(10)
$$ 
with $gcd(r,t)=1$.  If $gcd\left({{(4k+1)^{4m+2}-1}\over {8k}},
{{q_i^{2\alpha_i+1}-1}\over {q_i-1}}\right)=1$ for all $i$, 
the relation (10) requires ${{(4k+1)^{4m+2}-1}\over {8k}}\biggl\vert t^2$
or equivalently ${{(4k+1)^{4m+2}-1}\over {8k}}=\sigma_\ell \tau_\ell^2$ where
$\sigma_\ell \tau_\ell\vert t$.  The substitution $t=\sigma_\ell \tau_\ell u$ 
gives 
$$\eqalign{(4k+1)~\prod_{i=1}^\ell~{{q_i^{2\alpha_i+1}}\over {q_i-1}} \cdot
(\sigma_\ell \tau_\ell u)^2~&=~\sigma_\ell \tau_\ell^2 \cdot r^2
\cr
(4k+1)~\prod_{i=1}^\ell~{{q_i^{2\alpha_i+1}-1}\over {q_i-1}}\sigma_l u^2
~&=~r^2
\cr}
\eqno(11)
$$
which, in turn, requires that $(4k+1)~\prod_{i=1}^\ell~{{q_i^{2\alpha_i+1}-1}
\over {q_i-1}}~=~\sigma_l~v^2$ and $r=\sigma_l v u$, so that 
$\sigma_l u\vert r$ and $\sigma_l u\vert t$, contrary to the original 
assumption that $r$ and $t$ are relatively prime unless $\sigma_\ell=u=1$.   
The rationality condition reduces to the existence of solutions to the equation
$${{x^n-1}\over {x-1}}~=~2y^2~~~~~~~~~x\equiv 1~(mod~4),~n \equiv 2~(mod~4)
\eqno(12)
$$
This relation is equivalent to the two conditions
${{x^{2m+1}-1}\over {x-1}}=y_1^2,~{{x^{2m+1}+1}\over 2}=y_2^2,~
y=y_1y_2,~(y_1,y_2)=1$ since $gcd(x^{2m+1}-1,x^{2m+1}+1)=2$.  It can be 
verified that there are no integer solutions to these simultaneous
Diophantine equations, implying that when ${{(4k+1)^{4m+2}-1}\over {8k}}$ 
satisfies the gcd condition given above, the square root 
of 
\hfil\break
$\left[(4k+1)\prod_{i=1}^\ell
U_{2\alpha_i+1}(q_i+1,q_i){1\over 2}U_{4m+2}(4k+2,4k+1)\right]$ is not 
a rational number and 
\hfil\break
there is no odd perfect number of the form with 
this constraint on the pair $(4k+1,4m+2)$.
\vfill\eject
\noindent
{\bf 4. Prime power divisors of Lucas sequences and Catalan's conjecture}
 
The number of distinct prime divisors of ${{q^n-1}\over {q-1}}$ is bounded 
below by $\tau(n)-1$ if $q>2$, where $\tau(n)$ is the number of natural 
divisors of $n$ [1][15]. 
The characteristics of these prime divisors can be deduced from the
properties of Lucas sequences.  Since the repunits 
${{q_i^{2\alpha_i+1}-1}\over {q_i-1}}$ have only odd prime divisors, the 
proofs in the following sections will have general validity, circumventing 
any exceptions corresponding to prime $q=2$. 

For a primary recurrence relation, defined by the initial values $U_0~=~0$ and
$U_1~=~1$, denoting the least positive integer $n$ such that 
$U_n(a,b)~\equiv~0~(mod~p)$, the rank of apparition, by $\alpha(a,b,p)$,
it is known that $\alpha(x+1,x,p)=ord_p(x)$ [16]. 

The extent to which the arguments $a$ and $b$ determine the divisibility 
of $U_n(a,b)$ [17][18] can be summarized as follows:

Let $p$ be an odd prime.
\hfil\break
If $p\vert a$, $p\vert b$, then $p\vert U_n(a,b)$ for all $n > 1$.
\hfil\break
If $p\nmid a$, $p\vert b$, then either $p\vert U_n(a,b),~n\ge 1$ or
$p\nmid U_n(a,b)$ for any $n\ge 1$.
\hfil\break   
If $p\vert a$ and $p\nmid b$, then $p\vert U_n(a,b)$ for all even $n$ 
or all odd $n$ or $p\nmid U_n(a,b)$ for any $n\ge 1$.
\hfil\break
If $p\nmid a$, $p\nmid b$, $p\vert D=a^2-4b$, then $p\vert U_n(a,b)$
when $p\vert n$.
\hfil\break
If $p\nmid abD$, then $p\vert U_{p-\left({D\over p}\right)}(a,b)$.

For the Lucas sequence $U_n(q+1,q)$, there is no prime  
which divides both $q$ and $q+1$, and since only $q$ is a divisor of the
second parameter, there are no prime divisors of $U_n(q+1,q)$ from this
category because ${{q^n-1}\over {q-1}}\equiv 1~(mod~q)$.  If $p\vert (q+1)$,
then ${{q^n-1}\over {q-1}}\equiv {{1-(-1)^n}\over 2}\equiv 0~(mod~p)$ when
$n$ is even.  However, $p\nmid {{q^n-1}\over {q-1}}$ with $n$ odd, and
therefore, prime divisors from this class are not relevant for the study 
of the product of repunits with odd exponents.

When $a=q+1$ and $b=q$, $D=(q-1)^2$ and if $p\vert (q-1)$, then
$p\vert~U_n(q+1,q)$ when $p\vert n$.  However, 
$p^2\nmid {{q^p-1}\over {q-1}}$, and under this condition,
$p^2\nmid {{q^n-1}\over {q-1}}$ unless $n={\scriptstyle C} p^2$.
More generally, denoting the power of $p$ which exactly divides $a$ by 
$p^{v_p(a)}$, it can be deduced that $v_p\left({{q^n-1}\over {q-1}}
\right)=v_p(n)$ if $p\vert (q-1)$ and $\alpha(q+1,q,p)=p$ [9][15].  

From the last property, it follows that $\alpha(a,b,p)\vert \left(p-
\left(D\over p\right)\right)$
When $p\nmid (q-1)$, $\left({D\over p}\right)=1$ and $\alpha(q+1,q,p)
\vert (p-1)$. If ${p^2\not\bigg\vert}~{{q^{p-1}-1}\over {q-1}}$, then 
$\alpha(q+1,q,p^2)=p~\alpha(q+1,q,p)$ so that $\alpha(q+1,q,p^2)
\vert p(p-1)$.  If $p^2\bigg\vert {{q^{p-1}-1}\over {q-1}}$, 
$\alpha(q+1,q,p^2)=\alpha(q+1,q,p)\vert p-1$ [19][20].  Thus a repunit with
primitive divisor $p$ is also divisible by $p^2$ if $Q_q\equiv 0~(mod~p)$ 
where $Q_a={{a^{p-1}-1}\over p}$ is the Fermat quotient.  
  
Since $q^n-1=\prod_{d\vert n}\Phi_d (q)$ where $\Phi_n(q)$ is the
$n^{th}$ cyclotomic polynomial, it can be shown that the largest arithmetic
primitive factor [21]-[23]
of $q^n-1$ when $q\ge 2$ and $n\ge 3$ is
$$\eqalign{\Phi_n(q)&~~~~~~~~~~if~\Phi_n(q)~and~n~are~relatively~prime
           \cr
           {{\Phi_n(q)}\over p}&~~~~~~~~if~a~common~prime~factor~p~of
                                       ~\Phi_n(q)~and~n~exists
                                                       \cr}
\eqno(13)
$$
In the latter case, if $n= p^f p^{\prime f^\prime}
p^{\prime\prime f^{\prime\prime}}~...$ is the prime factorization
of $n$, then $\Phi_n(q)$ is divisible by $p$ if and only if
$e={n\over {p^f}}=ord_p(q)$ when $p\nmid (q-1)$, and moreover, 
$p\parallel \Phi_{ep^f}(q)$ when $f > 0$ [1].     

Division by $q-1$ does not alter the arithmetic primitive factor, since it  
is the product of the primitive divisors of $q^n-1$, which are also
the primitive divisors of ${{q^n-1}\over {q-1}}$.  For all primitive 
divisors, $p^\prime\nmid (q-1)$, so that $(p^\prime)^h\bigg\vert 
{{q^n-1}\over {q-1}}$ if $(p^\prime)^h\vert q^n-1$ and the arithmetic 
primitive factor again would include $(p^\prime)^h$.   The imprimitive 
divisors would be similarly unaffected because the form of the index 
$n=ep^f$ prevents $q-1$ from being a divisor of $\Phi_n(q)$ when 
$p\nmid (q-1)$.  If $p\vert (q-1)$, the rank of apparition for the 
Lucas sequence $\{U_n(q+1,q)\}$ is $p$, so that it is consistent to set 
$n=p^{f+1}$.  Then, $p\parallel \Phi_{p^{f+1}}(q)$ and the arithmetic 
primitive factor is ${{\Phi_{p^{f+1}}(q)}\over p}$.

If $(q_i-1)\nmid \Phi_{n_i}(q_i)$, the product of the arithmetic 
primitive factors of each repunit ${{q_i^{n_i}~-~1}\over {q_i-1}}$ and
${{(4k+1)^{4m+2}-1}\over {4k}}$ in the expression (5) is
$${{\Phi_{n_1}(q_1)}\over {p_1}}~{{\Phi_{n_2}(q_2)}\over {p_2}}~...~
{{\Phi_{n_\ell}(q_\ell)}\over {p_\ell}}~\times~\left[{{\Phi_{4m+2}(4k+1)}
\over {p_{\ell+1}}}\right]
\eqno(14)
$$ 
where the indices are odd numbers $n_i=2\alpha_i+1$, $p_i,~i=1,...,l$, 
represents the common factor of $n_i$ and $\Phi_{n_i}(q_i)$, and
$p_{\ell+1}$ is a common factor of $4m+2$ and $\Phi_{4m+2}(4k+1)$.
Division of $\Phi_{n_i}(q_i)$ by the prime $p_i$ is necessary only when
$gcd(n_i,\Phi_{n_i}(q_i))\ne 1$, and $p_i=P\left({{n_i}\over {gcd(3,n_i)}}
\right)$, where $P(n)$ represents the largest prime factor of $n$ [17]
[24]-[27]. 

\vskip 5pt
\noindent
{\bf Theorem 1.}  The arithmetic primitive factors of the repunits with
different prime bases could be equal only if the exponents are different,
with possible exceptions being determined by the solutions to the equation
${{q_j^n-1}\over {q_i^n-1}}=p,~q_i\ne q_j$ with $q_i, q_j$ and $p$ prime.
\vfill\eject
\noindent
{\bf Proof.}  Consider the following four cases:
\vskip 2pt
\noindent
I. The arithmetic primitive factors of $q_i^{n_i}-1$ and $q_j^{n_j}-1$ are 
$\Phi_{n_i}(q_i)$ and $\Phi_{n_j}(q_j)$.
\vskip 2pt
\noindent
Since $\Phi_n(x)$ is a strictly increasing function for $x\ge 1$ [28],
$\Phi_n(q_j)>\Phi_n(q_i)$ when $q_j$ is the larger prime, and
equality of $\Phi_{n_i}(q_i)$ and $\Phi_{n_j}(q_j)$ could only be
achieved, if at all feasible, when $n_i\ne n_j$.
\vskip 2pt
\noindent
II. The arithmetic primitive factors of $q_i^{n_i}-1$ and $q_j^{n_j}-1$ are
$\Phi_{n_i}(q_i)$ and  ${{\Phi_{n_j}(q_j)}\over {p_j}}$.
\vskip 2pt
\noindent
Comparing $\Phi_n(q_i)$ and ${{\Phi_n(q_j)}\over p}$, $p=p_j$ is a common 
factor of $n$ and $\Phi_n(q_j)$ but it does not divide $\Phi_n(q_i)$.  It 
follows that the relation $\Phi_n(q_i)={{\Phi_n(q_j)}\over p}$ could only 
hold if $p\parallel \Phi_n(q_j)$.  The prime decomposition of $e$ as
$\rho_1^{r_1}...\rho_s^{r_s},~gcd(\rho_t,p)=1,~t=1,...,s$, leads to the 
following expressions for $\Phi_n(q_i)$ and $\Phi_n(q_j)$, 

$$\eqalign{\Phi_n(q_i)&=\Phi_{ep^f}(q_i)={{\Phi_e(q_i^{p^f})}\over
{\Phi_e(q_i^{p^{f-1}})}}
\cr
&={{\prod_{{k~even}\atop {k\ge 0}}\prod_{{t_k>...>t_1}\atop {t_k\le s}}
\left[q_i^{{ep^f}\over {\rho_{t_1}...\rho_{t_k}}} -1\right]}
\over
{\prod_{{{\tilde k}~odd}\atop {{\tilde k}\ge 1}}\prod_{{t_{\tilde k}>...
>t_1}\atop {t_{\tilde k}\le s}}
~\left[q_i^{{ep^f}\over {\rho_{t_1}...\rho_{t_{\tilde k}}}}-1\right]}}
~\cdot~
{{\prod_{{{\tilde k}~odd}\atop {{\tilde k}\ge 1}}\prod_{{t_{\tilde k}>...>t_1}
\atop {t_{\tilde k}\le s}}
~\left[q_i^{{ep^{f-1}}\over {\rho_{t_1}...\rho_{t_{\tilde k}}}}-1\right]}
\over
{\prod_{{k~even}\atop {k\ge 1}}\prod_{{t_k>...>t_1}\atop {t_k\le s}}
\left[q_i^{{ep^{f-1}}\over {\rho_{t_1}...\rho_{t_k}}}-1\right]}}
\cr}
$$
$$\eqalign{\Phi_n(q_j)&=\Phi_{ep^f}(q_j)={{\Phi_e(q_j^{p^f})}\over
{\Phi_e(q_j^{p^{f-1}})}}
\cr
&={{\prod_{{k~even}\atop {k\ge 0}}
\prod_{{t_k>...>t_1}\atop {t_k\le s}}
\left[q_j^{{ep^f}\over {\rho_{t_1}...\rho_{t_k}}}-1\right]}
\over
{\prod_{{{\tilde k}~odd}\atop {{\tilde k}\ge 1}}\prod_{{t_{\tilde k}>...>t_1}
\atop {t_{\tilde k}\le s}}
~\left[q_j^{{ep^f}\over {\rho_{t_1}...\rho_{t_{\tilde k}}}}-1\right]}}
~\cdot~{{\prod_{{{\tilde k}~odd}\atop {{\tilde k}\ge 1}}
\prod_{{t_{\tilde k}>...>t_1}\atop {t_{\tilde k}\le s}}
~\left[q_j^{{ep^{f-1}}\over {\rho_{t_1}...\rho_{t_{\tilde k}}}}-1\right]}
\over
{\prod_{{k~even}\atop {k\ge 1}}\prod_{{t_k>...>t_1}\atop {t_k\le s}}
\left[q_j^{{ep^{f-1}}\over {\rho_{t_1}...\rho_{t_k}}}-1\right]}}
\cr}
\eqno(15)
$$

Since $e=ord_p(q_j)$, it follows that $p\vert (q_j^e-1)$, and if 
$q_j^e=1+pk_j~(mod~p)$, then $(q_j^e)^{p^f}=(1+pk_j)^{p^f}\equiv 1+p^fpk_j
\equiv 1~(mod~p^{f+1})$.  Thus, $p^{f+1}\vert (q_j^{ep^f}-1)$ and
$p^f\vert (q_j^{ep^{f-1}}-1)$, while ${p\not\bigg\vert}~
\left(q_j^{{e\over {\rho_{t_1}...\rho_{t_k}}}p^f}-1\right)$. 
Let $H(f)\ge f+1$ denote the exponent such that 
\vfill\eject
\noindent
$p^{H(f)}\parallel (q_j^{ep^f}-1)$.  Since $q_j^{ep^{f-1}}
\equiv 1~(mod~p^{H(f-1)})$, $q_j^{ep^f}= (1+p^{H(f-1)}k^\prime_j)^p
\equiv 1+p\cdot p^{H(f-1)}k^\prime_j\equiv 1~(mod~p^{H(f-1)+1})$.  
Consequently, $H(f)-H(f-1)=1$, which is consistent with $\Phi_n(q_j)$ being 
exactly divisible by $p$.  

Although $\Phi_n(q_i)$ and ${{\Phi_n(q_j)}\over p}$ are not divisible by 
$p$, consider a primitive prime factor $p^\prime$ of $\Phi_n(q_i)$.  It must 
divide some factor $q_i^{n\over {\rho_{t_1}...\rho_{t_k}}}-1$ in the 
expression for $\Phi_n(q_i)$, and thus, it will also divide $q_i^{n\over 
{\rho_{t_1}...\rho_{t_\ell}}}-1,~\ell<k$.  Since the exponent of
$q_i^{n\over {\rho_{t_1}...\rho_{t_\ell}}}-1$ in $\Phi_n(q_i)$ is
$(-1)^l$, there will be $2^{k-1}$ factors in the numerator and $2^{k-1}$
factors in the denominator divisible by $p^\prime$.  When $k\ge 1$, the
factors of $p^\prime$ exactly cancel because each term $q_i^{n\over
{\rho_{t_1}...\rho_{t_\ell}}}-1$ is divisible by the same power of $p^\prime$.
The exception occurs when $p^\prime\vert q_i^n-1$ only; if
$p_a^{\prime f_a}\parallel q_i^n-1$, then $p_a^{\prime f_a}\parallel 
\Phi_n(q_i)$ [29].  Equivalence of $\Phi_n(q_i)$ and ${{\Phi_n(q_j)}\over p}$
requires that the prime power divisors of these quantities are equal,
so that $p_a^{\prime f_a}\bigg\vert\bigg\vert {{\Phi_n(q_j)}\over p}$ for all
primes $\{p^\prime_a\}$.  However, if $p_a^{\prime f_a}\parallel q_j^n-1$,
then $q_i^n-1$ and $q_j^n-1$ have the same primitive prime power divisors.
The imprimitive prime divisor $p$ which divides $q_j^n-1$ might also
divide $q_i^n-1$, although overall cancellation of $p$ in $\Phi_n(q_i)$ 
requires that $p^r\bigg\vert q_i^{n\over {\rho_{t_1}...\rho_{t_k}}}-1$ for 
some $k\ge 1$ and $p^{r-1}\bigg\vert 
q_i^{n\over {p\rho_{t_1}...\rho_{t_k}}}-1$.  When $p^r\parallel q_i^n-1$
and ${{q_j^n-1}\over {q_i^n-1}}=p^{H(f)-r}$ 
$$q_i^n-1=\kappa u_1^{H(f)-r}~~~~~q_j^n-1=\kappa u_2^{H(f)-r}~~~~~
{{u_2}\over {u_1}}=p
\eqno(16)
$$

Integer solutions of $w=y^m,~y\ge 2,~m\ge 2$ can be written as 
$w=x^n,~x\ge 2$ with $m\vert n$.  Since $y\vert x^n$, $y\nmid (x^n-1)$ 
because $y\ge 2$.   The nearest integers to $x^n$ having a similar form,
$\{(x-1)^n, (x+1)^n, (x+1)^{n-1}, (x-1)^{n+1}\}$ do not provide a
counterexample to the conclusion since none of them are divisible by $y$.
Furthermore, $x^n-(x-1)^n~>~1$, $(x+1)^n-x^n~>~1$, $\vert (x+1)^{n-1}
-x^n\vert >1,~x\ge 2,n\ge 4; x\ge 3, n\ge 3$ and $\vert x^n-(x-1)^{n+1}
\vert~>~1,~x\ge 2, n\ge 3$ so that none of these integers will have the 
form $y^m\pm 1$.  The exception occuring when $x=y=2,~m=n=3$ is the statement 
of Catalan's conjecture, that $(X,Y,U,V)=(3,2,2,3)$ is the only 
integer solution of $X^U-Y^V=1$.  Thus, if $\kappa=1$, any non-trivial 
solution to equation (16) is constrained by the condition $H(f)-r=1$, which 
implies that ${{q_j^n-1}\over {q_i^n-1}}=p$.  Since the odd primes $q_i,q_j$ 
and the exponent $n$ in the prime decomposition of $N$ must be greater 
than or equal to 3, this restriction is consistent with Catalan's conjecture.  

When $\kappa\ne 1$, it may be noted that for $q_i,~q_j \gg 1$, 
${{q_j^n-1}\over {q_i^n-1}}\simeq \left({{q_j}\over {q_i}}\right)^n\ne p^h$.  
Exceptional solutions to equation (16) occur, for example, when $h=1$; 
they include $\{(q_i,q_j;n;p)= (3,5;2;3),~(5,7;2;2),
~(5,11;2;5),~(5,13;2;7),~(11,19;2;3),~(7,23;2;11),~(11,29;2;7),$
\hfil\break
$~(29,41;2;2)\}$.  
Since $q_i\ne q_j$, with the exception of the non-trivial solutions to 
equation (16), it would be necessary to set $n_i\ne n_j$ to obtain equality 
between $\Phi_{n_i}(q_i)$ and ${{\Phi_{n_j}(q_j)}\over p}$. 
\vskip 2pt
\noindent
III.  The arithmetic primitive factors of $q_i^{n_i}-1$ and $q_j^{n_j}-1$
are ${{\Phi_{n_i}(q_i)}\over {p_i}}$ and $\Phi_{n_j}(q_j)$.
\vskip 2pt
\noindent
The proof of the necessity of $n_i\ne n_j$ for any equality between the
arithmetic primitive factors is similar to that given in Case II with the
roles of $i$ and $j$ interchanged.
\vskip 2pt
\noindent
IV.  The arithmetic primitive factors of $q_i^{n_i}-1$ and $q_j^{n_j}-1$
are ${{\Phi_{n_i}(q_i)}\over {p_i}}$ and ${{\Phi_{n_j}(q_j)}\over {p_j}}$.
\vskip 2pt
\noindent
Since $p_i=gcd(n_i,\Phi_{n_i}(q_i))$ and $p_j=gcd(n_j,\Phi_{n_j}(q_j))$, 
$\Phi_{n_i}(q_i)$ and $\Phi_{n_j}(q_j)$ share a common factor if $n_i=n_j$.  
Thus, the primes $p_i$ and $p_j$ must be equal, and a comparison can be made
between ${{\Phi_n(q_i)}\over p}$ and ${{\Phi_n(q_j)}\over p}$.  Again,
by the monotonicity of $\Phi_n(x)$, it follows that these quantities are
not equal when $q_i$ and $q_j$ are different primes.  Equality of the
arithmetic prime factors could only occur if $n_i\ne n_j$.\qed

\vskip 10pt
\noindent
{\bf 5. The exponent of prime divisors of repunit factors in the
rationality condition}

Since all primitive divisors of $U_n(a,b)$ have the form $p=nk+1$, it
follows that $p\bigg\vert {{q^{{(p-1)}\over {\iota(p)}}-1}\over {q-1}}$.  If
$\iota(p)$ is odd, where $\iota(p)$ is the residue index, the exponent 
${{p-1}\over {\iota(p)}}$ will be even for all odd primes $p$, whereas if 
$\iota(p)$ is even, the exponent ${{p-1}\over {\iota(p)}}$ may be 
even or odd.  Given that $p\vert U_{2\alpha_i+1}(q_i+1,q_i)$, 
$\iota(p)$ is even and $p\bigg\vert {{q_i^{{(p-1)}\over 2}-1}
\over {q_i-1}}$ implying $q_i^{{p-1}\over 2}\equiv 1~(mod~p)$ and 
$\left({{q_i}\over p}\right)=1$.  Moreover, if $\left({{q_i}\over p}\right)
=\left({{q_j}\over p}\right)=1$, $\left({{q_iq_j}\over p}\right)=1$ 
implying that $p\vert {(q_iq_j)}^{{p-1}\over 2}-1$.  Thus, the 
Fermat quotient is $Q_{q_iq_j}={{{(q_iq_j)}^{{p-1}\over 2}-1}\over p}
\left({(q_iq_j)}^{{p-1}\over 2}+1\right)={\cal N}_{q_iq_j}
({\cal N}_{q_iq_j}p+2)$ where ${\cal N}_q$ can be defined to be 
${{q^{{(p-1)}\over 2}-1}\over p}$.  By the logarithmic rule for Fermat 
quotients, $Q_{qq^\prime}\equiv Q_q+Q_{q^\prime}~(mod~p)$ [30], so that 
${\cal N}_{q_iq_j}\equiv {\cal N}_{q_i}+{\cal N}_{q_j}~(mod~p)$.

Recalling that $\alpha(q_i+1,q_i,p^2)\ne\alpha(q_i+1,q_i,p)$ only when 
$p^2\not\bigg\vert~{{q_i^{p-1}-1}\over {q_i-1}}$, it is sufficient to
prove that the Fermat quotient $Q_{q_i}\ne 0~(mod~p)$ to show that the
$p^2$ is not a divisor of the repunit ${{q_i^{2\alpha_i+1}-1}\over {q_i-1}}$.
It has been established that $q^{p-1}-1\equiv p\bigg({\mu_1}+{{\mu_2}\over 2}
+...+{{\mu_{p-1}}\over {p-1}}\bigg)$
\hfil\break
$~(mod~p^2)$, where $\mu_i\equiv \left[{{-i}\over p}\right]
~(mod~q)$ is the minimum positive integer congruent to $(-i\cdot p^{-1})~
(mod~q)$ [31][32].  Since $\mu_i\ne 0$ in general, except when $i=q$, it
follows that $q^{p-1}-1\ne 0~(mod~p^2)$ except for $p-1$ values of 
$q$ between 1 and $p^2-1$.     

By Hensel's lemma [33][34], each of the integers between 1 and 
$p-1$, which satisfy $x^{p-1}-1\equiv 0~(mod~p)$ generate
the $p-1$ solutions to the congruence equation 
$$(x^{\prime})^{p-1}-1\equiv 0~(mod~p^2)
\eqno(17)
$$
through the formula
$$x^\prime~=~x+\left({{-g_1(x)p}\over {(p-1)q^{p-2}}}\right)~(mod~p^2)
\eqno(18)
$$
with $x^{p-1}-1\equiv g_1(x)p~(mod~p^2)$.  

Since $\varphi(p^2)=p(p-1)$, a set of $p-1$ solutions to equation (17) 
can also be labelled as $c^p~(mod~p^2),~1\le c\le p-1$, 
since $(c^p)^{p-1}=c^{p(p-1)}=c^{\varphi(p^2)}\equiv 1~(mod~p^2)$.  Each 
power $c^p$ is different, because $c_1^p\equiv c_2^p~(mod~p^2)$ 
implies $c_1=c_2$ since $p^2\nmid(c_3^p-1)$ for any $c_3$ between 
1 and $p-1$.

Theorems concerning the Fermat quotient ${{q^r-1}\over {q-1}}$ can be 
extended to quotients of the type ${{q^{nr}-1}\over {q^n-1}}$.  It has been
proven, for example, that $p\bigg\vert\bigg\vert 
{{q^{nr}- 1}\over {q^n-1}}$,  $p\nmid r,~p\nmid q^n-1$, 
then $Q_q={{q^{p-1}-1}\over p}\not\equiv 0~(mod~p)$ [35], and more generally, 
if $p^h~\bigg\vert\bigg\vert~{{q^{nr}-1}\over {q^n-1}}$, $p\nmid r$, 
$p~\nmid q^n-1$, then $q^{p-1}\nmid 1 (mod~p^{h+1})$. 

When $p\vert (q^r-1)$, the following lemma
is obtained. 
\vskip 10pt
\noindent
{\bf Lemma.} For any prime $p$ which is a primitive divisor of
$U_{2\alpha_i+1}(q_i+1,q_i)$, $p\not\bigg\vert {{q_i^{p-1}-1}\over 
{q_i^{(2\alpha_i+1)}-1}}$, and if $p^h\parallel U_{2\alpha_i+1}(q_i+1,q_i)$, 
then $p^h\bigg\vert\bigg\vert {{q_i^{2\alpha_i+1}
-1}\over {q_i^{{(2\alpha_i+1)}\over s}-1}}$ for any non-trivial divisor 
$s$ of $2\alpha_i+1$.
\vskip 5pt
\noindent
{\bf Proof.} 
Defining the residue index $\iota_i(p)$ by $p-1=(2\alpha_i+1)\iota_i(p)$, then
$$p~\bigg\vert~{{q_i^{p-1}-1}\over {q_i-1}}=
\left[{{q^{(2\alpha_i+1)\iota_i(p)}-1}\over {q_i^{(2\alpha_i+1)}-1}}
\right]\cdot\left[{{q^{(2\alpha_i+1)}-1}\over {q_i-1}}\right]
\eqno(19)
$$
Suppose that $p \vert {{q_i^{(2\alpha_i+1)\iota_i(p)}-1} \over
{q_i^{2\alpha_i+1}-1} }$.  Then, by equation (19), 
$p^2\bigg\vert {{q_i^{p-1}-1}\over {q-1}}$.  By a theorem on congruences, if 
$q^e\equiv 1~(mod~p)$, where $e\vert (p-1)$ and 
$q^{p-1}\equiv 1~(mod~p^2)$, then $q^e\equiv 1~
(mod~p^2)$ [28], so that $p^2\bigg\vert {{q_i^{2\alpha_i+1}-1}\over {q_i-1}}$.
Consequently, $p^3\bigg\vert {{q_i^{p-1}-1}\over {q_i-1}}$.  The 
theorem on congruences can be extended to larger prime powers:
$q^e\equiv 1~(mod~p^n)$ and $q^{p-1}\equiv 1~(mod~p^{n+1})$, then
$q^e\equiv 1~(mod~p^{n+1})$.  From the first congruence relation,
$q^e=1+k^\prime p^n$ for some integer $k^\prime$.  Raising this quantity to
the power ${{p-1}\over e}$, it follows that
$$1~\equiv~q^{p-1}=(q^e)^{{p-1}\over e}=(1+k^\prime p^n)^{{p-1}\over e}
\equiv 1+k^\prime p^n {{p-1}\over e}~(mod~p^{n+1})
\eqno(20)
$$
Since ${{p-1}\over e}<p$, the integer $k^\prime$ must be a multiple of 
$p$.  Thus, $q^e=1+k^{\prime\prime}p^{n+1}\equiv 1~(mod~p^{n+1})$.
By the generalized congruence theorem, $p^3\bigg\vert 
{{q_i^{2\alpha_i+1}-1}\over {q_i-1}}$ and equation (20) 
in turn implies that $p^4\bigg\vert {{q_i^{p-1}-1}\over {q_i-1}}$.  
Since this process can be continued indefinitely to arbitrarily high 
powers of the prime $p$, a contradiction is obtained once the maximum 
exponent is greater than $h$, where $p^h\bigg\vert {{q_i^{p-1}-1}\over 
{q_i-1}}$.  Therefore,  $p\not\bigg\vert~{{q_i^{(2\alpha_i+1)\iota_i(p)}-1}
\over {q_i^{2\alpha_i+1}-1}}$. 

Similarly,
$${{q_i^{2\alpha_i+1}-1}\over {q_i-1}}=\left[{{q_i^{2\alpha_i+1}-1}\over
{q_i^{{2\alpha_i+1}\over s}}-1}\right]\cdot
\left[{{q_i^{{2\alpha_i+1}\over s}-1}\over {q_i-1}}\right]
\eqno(21)
$$ 
If $s$ is a non-trivial divisor of $2\alpha_i+1$, then 
$p\not\bigg\vert~{{q_i^{{2\alpha_i+1}\over s}-1}\over {q_i-1}}$, because it
is a primitive divisor of $U_{2\alpha_i+1}(q_i+1,q_i)$.  Given that
$p^h\parallel U_{2\alpha_i+1}(q_i+1,q_i)$, by equation (21), 
$p^h\bigg\vert\bigg\vert {{q_i^{(2\alpha_i+1)}-1}\over 
{q_i^{{(2\alpha_i+1)}\over s}-1}}$. \qed

Imprimitive prime divisors of $U_n(a,b)$ are characterized by the 
property that $p\vert U_d(a,b)$ for some $d\vert n$.  The exponent of the 
imprimitive prime power divisor exactly dividing ${{q^n-1}\over {q-1}}$ 
can be determined by a further lemma: if $p^h\bigg\vert {{q^n-1}\over {q-1}}$,
then either $gcd(n,p-1)=1$, $q\equiv 1~(mod~p)$, $p^h\vert n~(mod~p)$ or
$e=gcd(n,p-1)>1$, $p^k\vert \Phi_e(q)$, $p^{h-k}\parallel n$ [15].
Since $v_p(\Phi_e(q))=v_p(q^e-1)$ if $p\nmid q-1$, the 
general formula [12] for the exponent of a prime divisor of a repunit
is
$$v_p\left({{q^n-1}\over {q-1}}\right)~=~\Bigg\{{{ {{v_p(q^e-1)+v_p(n)~~~~~~
~~~e=ord_p(q)\vert n,~e>1}\atop {v_p(n)~~~~~~~~~~~~~~~~~~~~~p\vert q-1}}}
\atop {0~~~~~~~~~~~~~~~~~~~~~~\textstyle{otherwise}}}
\eqno(22)$$
The exponent also can be deduced from the congruence properties of 
$q$-numbers $[n]={{q^n-1}\over {q-1}}$ and $q$-binomial coefficients [36],
as it equals $s=\epsilon_0h+\epsilon_1+...+\epsilon_{k-1}$ where $p^h\parallel
q^e-1$ and
$$\eqalign{n-1&=a_0+e(a_1+a_2p+...+a_kp^{k-1})
\cr
n&=b_0+e(b_1+b_2p+...+b_kp^{k-1})
\cr
~a_0\le e-1,&~a_i\le p-1,~i=1,...,k-1 
\cr
~b_0\le e-1,&~b_i\le p-1,~i=1,...,k-1
\cr
a_0+1~&=~\epsilon_0 e+b_0
\cr
\epsilon_0+a_1~&=~\epsilon_1 p +b_1
\cr
&~...
\cr
\epsilon_{k-2}+a_{k-1}~&=~\epsilon_{k-1}p+b_{k-1}
\cr
\epsilon_{k-1}+a_k~&=~b_k
\cr}
\eqno(23)
$$
with $\epsilon_i$ equal to 0 or 1, which is consistent with equation (21) 
because $\epsilon_0=1$ and $v_p(n)=\epsilon_1+...+\epsilon_{k-1}$.

Specializing to the case of $h=2$, it follows that if the quotient
${{q^n-1}\over {q-1}}$ is exactly divisible by $p^2$, then
\vskip 2pt
\noindent
(i) $gcd(n,p-1)=1$, $p\vert (q-1)$ or $p\nmid q^n-1$, $p^2\parallel n$ 
\hfil\break
(ii) $p\parallel \Phi_e(q)$, where $e=\alpha(q+1,q,p)$ is the rank of
apparition of $p$, $p\parallel n$ 
\hfil\break
(iii) $p^2\parallel \Phi_e(q)$, $p\nmid n$   
\vskip 2pt
\noindent
and the only indices $n_i$ which allow for exact divisibility of 
${{q_i^{n_i}-1}\over {q_i-1}}$ by $p^2$ are $n_i=\mu p^2$, when 
$p\vert (q_i-1)$ or $e_i\nmid n_i$, $n_i=\mu e_i p$ when 
$p\parallel \Phi_{e_i}(q_i)$ and $n_i=\mu e_i$ when $p^2\parallel 
\Phi_{e_i}(q_i)$.  Since $n_i$ is odd,
the three categories can be defined by the conditions: (i) $n_i=\mu p^2$,
(ii) $n_i=\mu e_i p$, $p$ is a primitive divisor of ${{q_i^{e_i}-1}\over 
{q_i-1}}$, $Q_{q_i}\not\equiv 0~(mod~p)$ 
(iii) $p$ is a primitive divisor of ${{q_i^{e_i}-1}\over {q_i-1}}$, 
$Q_{q_i}\equiv 0~(mod~p)$.
\vskip 10pt
\noindent
{\bf 6. A proof by the method of induction of the non-existence of a generic
set of} 
\hfil\break
\phantom{.....}{\bf primes satisfying the rationality condition}

The equation
$$a{{x^m-1}\over {x-1}}~=~b{{y^n-1}\over {y-1}}
\eqno(24)
$$
is known to have finitely many integer solutions for $m,n,x,y$, given
$a$ and $b$ such that $gcd(a,b)=1$, $a(y-1)\ne b(x-1)$, and 
$max(m,n,x,y)< C$ where $C$ is an effectively computable number depending
on $a,b$ and $F$ where $\vert x-y\vert < F {z\over {(log~z)^2(log~log~z)^3}}$ 
with $z=max(x,y)$ [37][38].   Using this relation to re-express 
${{q_i^{2\alpha_i+1}-1}\over {q_i-1}}$ in terms of 
${{(4k+1)^{4m+2}-1}\over {4k}}$, it can be established that there 
are unmatched primes in the product of the repunits (5) and 
that the square root of this expression is irrational for several different 
categories of prime divisors $\{q_i,~i=1,...,\ell;~4k+1\}$. 
\hfil\break\hfil\break
{\bf Theorem 2.}  The square-root expressions ${\sqrt {2(4k+1)}}
\left[{{q_1^{2\alpha_1+1}-1}\over {q_i-1}}...{{q_\ell^{2\alpha_\ell+1}-1}
\over {q_\ell-1}}\right]^{1\over 2}$
\hfil\break
$\cdot \left({{(4k+1)^{4m+2}-1}\over {4k}}\right)^{1\over 2}$ are not 
rational numbers for the following sets of primes $\{q_i,i=1,...,\ell;$
$~4k+1\}$ and exponents $2\alpha_i+1$:
\hfil\break\hfil\break
(i) For sets of primes with the number of elements given by
consecutive integers, $\{q_i, i=1,...\ell-1,~4k+1\}$ and 
$\{q_j^\prime,~j=1,...,\ell,~4k^\prime+1\}$, there cannot be odd integers of
the form $N_1=(4k+1)^{4m+1}q_1^{2\alpha_1}...q_{\ell-1}^{2\alpha_{\ell-1}}$
and $N_2=(4k^\prime+1)^{4m^\prime+1}(q_1^\prime)^{2\alpha_1^\prime}
...(q_\ell^\prime)^{2\alpha_\ell^\prime}$ such that both ${{\sigma(N_1)}\over
{N_1}}=2$ and ${{\sigma(N_2)}\over {N_2}}=2$.
\hfil\break\hfil\break
(ii) Setting $\alpha_j=\alpha_\ell$, extra prime divisors $p$ of the
repunits ${{q_j^{2\alpha_j+1}-1}\over {q_j-1}},~j<\ell$ and 
${{q_\ell^{2\alpha_\ell+1}-1}\over {q_\ell-1}}$, where $p\vert~(q_j-1)$ but 
$p\nmid (q_\ell-1)$, cannot be absorbed into the square factors if 
\hfil\break\hfil\break
$Q_{q_\ell}\not\equiv~0~(mod~p)$ or $p^{h_\ell^\prime}
\bigg\vert\bigg\vert {{q_\ell^{2\alpha_\ell+1}-1}\over {q_\ell-1}}$ with 
$h_\ell^\prime$ odd.  Similarly, if $p\nmid (q_j-1)$ but 
$p\vert~(q_\ell-1)$, then an 
odd power of $p$ divides the product of the two repunits if 
$Q_{q_j}\not\equiv 0~(mod~p)$ or $Q_{q_j}\equiv 0~(mod~p^{h_j^\prime-1})$,
$Q_{q_j}\not\equiv 0~(mod~p^{h_j^\prime})$, with $h_j^\prime$ odd, 
and $p$ remains an unmatched prime divisor.
\hfil\break\hfil\break
(iii) When $n_j=2\alpha_j+1$ is set equal to $n_\ell=2\alpha_\ell+1$, the 
primitive prime divisors of ${{q_j^{2\alpha_j+1}-1}\over {q_j-1}}$
and ${{q_\ell^{2\alpha_\ell+1}-1}\over {q_\ell-1}}$ cannot not be matched
to produce the square of a rational number if ${{q_\ell^{n_\ell}-1}
\over {q_j^{n_\ell}-1}}\ne {{y_2^2}\over {y_1^2}},~y_1,~y_2~\in~{\Bbb Z}$.
This property is valid, for example, when 
$q_\ell^{{n_\ell}\over 2}< gcd(q_j^{n_j}-1, q_\ell^{n_\ell}-1)$.
\hfil\break\hfil\break
(iv) Additional prime divisors are introduced when the exponents are adjusted,
so that, in general, there will be unmatched prime divisors in the products of 
the repunits ${{q_j^{2\alpha_j+1}-1}\over {q_j-1}}$,
\hfil\break
$~j<\ell$ and ${{q_\ell^{2\alpha_\ell+1}-1}\over {q_\ell-1}}$ when 
$\alpha_j\ne \alpha_\ell$.
\vfill\eject
\noindent
{\bf Proof.}
\hfil\break
Suppose $\{a_i\}$ and $\{b_i\}$ are defined by 
$$\eqalign{a_1{{q_1^{2\alpha_1+1}-1}\over {q_1-1}}~=&~b_1{{(4k+1)^{4m+2}
                                                           -1}\over {4k}}
\cr
a_2{{q_2^{2\alpha_2+1}-1}\over {q_2-1}}~=&~b_2{{(4k+1)^{4m+2}-1}\over {4k}}
\cr
\vdots&
\cr}
$$
$$
a_\ell{{{q_\ell}^{2\alpha_\ell+1}-1}\over {q_\ell-1}}~=~b_\ell
{{(4k+1)^{4m+2}-1}\over {4k}}
\eqno(25)
$$
Then 
$$\eqalign{{\sqrt {2(4k+1)}}&\left[{{q_1^{2\alpha_1+1}-1}\over {q_1-1}} 
{{q_2^{2\alpha_2+1}-1}\over {q_2-1}}~...~{{{q_\ell}^{2\alpha_\ell+1}-1}
\over {q_\ell-1}}{{(4k+1)^{4m+2}-1}\over {4k}}\right]^{1\over 2}
\cr
&~~~~~~~~~~~~~~=~
{\sqrt {2(4k+1)}}~{{(b_1b_2....b_\ell)^{1\over 2}}\over 
{(a_1a_2...a_\ell)^{1\over 2}}}\cdot\left({{(4k+1)^{4m+2}-1}\over {4k}}
                          \right)^{{(\ell+1)}\over 2}
\cr}
\eqno(26)
$$

If
$$a_{ij}{{q_i^{2\alpha_i+1}-1}\over {q_i-1}}~=~b_{ij}{{q_j^{2\alpha_j+1}-1}
\over {q_j-1}}
\eqno(27)
$$
define $\{a_{ij}\}$ and $\{b_{ij}\}$ with $gcd(a_{ij},b_{ij})=1$,
$${{b_1b_2b_3}\over {a_1a_2a_3}}~=~{{b_{13}}\over {a_{13}}}{{a_2}\over {b_2}}
\times \left({{b_2b_3}\over {a_2a_3}}\right)^2
\eqno(28)
$$
and
$${{b_1b_2...b_\ell}\over {a_1a_2...a_\ell}}~=~{{b_{1\ell}}\over {a_{1\ell}}}
{{a_2}\over {b_2}}...{{a_{\ell-1}}\over {b_{\ell-1}}}
\left({{b_2b_3...b_\ell}\over {a_2a_3...a_\ell}}\right)^2
\eqno(29)
$$
Since the fraction ${{b_1}\over {a_1}}$ can be expressed in terms of ${{b_2}
\over {a_2}}$
$${{b_1}\over {a_1}}~=~{{b_2}\over {a_2}}{{\rho_{12}^{r_{12}}}\over 
{\chi_{12}^{s_{12}}}}~=~{{b_2}\over {a_2}}{{\rho_{12}^{(r_{12})_0}}\over
{\chi_{12}^{(s_{12})_0}}}{{\rho_{12}^{(r_{12}-(r_{12})_0)}}\over
{\chi_{12}^{(s_{12}-(s_{12})_0)}}}
\eqno(30)
$$ 
where $\rho_{12}^{r_{12}},~\chi_{12}^{s_{12}}$ denote products of various 
powers of different primes, with $r_{12},~s_{12}$ representing the sets
of exponents, $(r_{12})_0,~(s_{12})_0$ labelling a collection of exponents 
consisting of 0 or 1 and  $\rho_{12},\chi_{12}$ being products of these 
primes with all of the exponents equal to 1.  The sets 
$(r_{12})_0,~(s_{12})_0$ are chosen so that 
$r_{12}-(r_{12})_0=2{\bar r}_{12},~s_{12}-(s_{12})_0=
2{\bar s}_{12}$ represent even exponents.  Since a similar relation
exists between ${{a_2}\over {b_2}}$ and ${{a_3}\over {b_3}}$,

$${{b_{13}}\over {a_{13}}}{{a_2}\over {b_2}}~=~\left[{{a_2}\over {b_2}}
{{\rho_{12}^{(r_{12})_0}}\over {\chi_{12}^{(s_{12})_0}}} 
{{\rho_{23}^{(r_{23})_0}}
\over {\chi_{23}^{(s_{23})_0}}}\right] \left({{\rho_{12}^{{\bar r}_{12}}}
\over {\chi_{12}^{{\bar s}_{12}}}}\right)^2 
\left({{\rho_{23}^{{\bar r}_{23}}}
\over {\chi_{23}^{{\bar s}_{23}}}}\right)^2
\eqno(31)
$$ 

If $(r_{12})_0=(s_{12})_0=(r_{23})_0=(s_{23})_0=\{0\}$, then ${{b_{13}}\over
{a_{13}}}{{a_2}\over {b_2}}\ne {{4k+1}\over 2}\cdot \square$ \footnote{*}
{$\square$ denotes the square of a rational number.} because rationality
of ${\sqrt{2(4k+1)}}\left[{{q_2^{2\alpha_2+1}-1}\over {q_2-1}}
{{(4k+1)^{4m+2}-1}\over {4k}}\right]^{1\over 2}$ would imply that
${{\sigma(N)}\over N}=2$ where $N=(4k+1)^{4m+1}q_2^{2\alpha_2}$,
contradicting the non-existence of odd perfect numbers with two prime
factors.

If $(r_{12})_0=(s_{12})_{0}=\{0\}$, the expression in brackets is not  
${{4k+1}\over 2}$ times the square of a rational number because 
$${{a_2}\over {b_2}} {{\rho_{23}^{(r_{23})_0}}\over {\chi_{23}^{(s_{23})_0}}}
= {{a_2}\over {b_2}} {{\rho_{23}^{r_{23}}}\over {\chi_{23}^{s_{23}}}}
\cdot \left({{\rho_{23}^{2{\bar r}_{23}}}\over {\chi_{23}^{2{\bar s}_{23}}}}
\right)^{-1}={{a_3}\over {b_3}}\left({{\chi_{23}^{{\bar s}_{23}}}\over
{\rho_{23}^{{\bar r}_{23}}}}\right)^2
\eqno(32)
$$
and ${{a_3}\over {b_3}}\ne {{4k+1}\over 2}\cdot \square$, since
${\sqrt{2(4k+1)}}\left[{{q_3^{2\alpha_3+1}-1}\over {q_3-1}}
{{(4k+1)^{4m+2}-1}\over {4k}}\right]^{1\over 2}$ is not rational and
there is no odd integer $N$ of the form $(4k+1)^{4m+1}q_3^{2\alpha_3}$
such that ${{\sigma(N)}\over N}=2$.  A similar conclusion holds 
when $(r_{23})_0=\{0\}$ and $(s_{23})_0=\{0\}$.

If both fractions ${{\rho_{12}^{(r_{12})_0}}\over 
{\chi_{12}^{(s_{12})_0}}}$, ${{\rho_{23}^{(s_{23})_0}}\over 
{\chi_{23}^{(s_{23})_0}}}$ are non-trivial, at least one of the pair
of exponents ($(r_{12})_0,~(s_{12})_0$), and at least one
of the pair of exponents ($(r_{23})_0,~(s_{23})_0)$, must equal one.
Under these conditions, the argument is not essentially changed when
all of the exponents are set equal to one, because replacement of the
prime factors in any of the coefficients $\rho_{12}$, $\chi_{12}$,
$\rho_{23}$ and $\chi_{23}$ by 1 only eliminates the presence of
these prime factors from the remainder of the proof.  The non-triviality
of both fractions, therefore, can be included by setting $(r_{12})_0
=(r_{23})_0=\{1\}$ and $(s_{12})_0=(s_{23})_0=\{1\}$.  The expression (32) 
then would be ${{4k+1}\over 2}$ times the square of a rational number if
$$\eqalign{a_2~=~(4k+1)\rho_{12}\cdot \rho_{23}\cdot {{p^2}\over 2}&~~~~~~~
                   b_2~=~\chi_{12}\cdot \chi_{23}\cdot q^2
\cr
&~~or
\cr
a_2~=~(4k+1)\chi_{12}\cdot \chi_{23}\cdot {{p^2}\over 2}&~~~~~~~
b_2~=~\rho_{12}\cdot\rho_{23}\cdot q^2
\cr}
\eqno(33)
$$
where $gcd(p,q)=1$.  If $a_2~=~(4k+1)\rho_{12}\rho_{23}{{p^2}\over 2}$ and
$b_2~=~\chi_{12}\chi_{23}q^2$,
$$\eqalign{{{a_3}\over {b_3}}~=~{{a_2}\over {b_2}}{{\rho_{23}}\over 
{\chi_{23}}}
~&=~{{(4k+1)\rho_{12}\rho_{23}^2p^2}\over {2\chi_{12}\chi_{23}^2q^2}}
\cr
{{a_3}\over {b_3}}{{2\chi_{12}}\over {(4k+1)\rho_{12}}}~&=~{{(\rho_{23}p)^2}
\over {(\chi_{23}q)^2}}
\cr}
\eqno(34)
$$
Since $gcd(a_3,b_3)=1$, the square-free factors can be separated in the 
fraction ${{a_3}\over {b_3}}={{{\hat a}_3}\over 
{{\hat b}_3}}\cdot {{{\hat p}^2}\over {{\hat q}^2}}$,
$${{2\chi_{12}}\over {(4k+1)\rho_{12}}}{{{\hat a}_3}\over {{\hat b}_3}}~
=~{{(\rho_{23}p{\hat q})^2}\over {(\chi_{23}q{\hat p})^2}}
\eqno(35)
$$
Since $a_3$ is even, and ${\hat a}_3$ is divisible by a single factor of 
$2$, $\chi_{12}=\rho_{23}{p\over 2}{\hat q}$, and similarly, because
${\hat b}_3$ is odd, $\rho_{12}={1\over {4k+1}}\chi_{23}q{\hat p}$.  Since 
${{\rho_{12}\rho_{23}}\over {\chi_{12}\chi_{23}}}= {2\over {4k+1}}
{{q{\hat p}}\over {p{\hat q}}}$, rationality of 
$\left[{2\over {4k+1}}{{b_{13}}\over {a_{13}}}
{{a_2}\over {b_2}}\right]^{1\over 2}$ also could be achieved by setting
$a_2=(4k+1)q {\hat p}{{p^{\prime 2}}\over 2}$ and 
$b_2= p {\hat q}q^{\prime 2}$.  Then
$$\eqalign{{{a_2}\over {b_2}}\cdot {{\rho_{23}}\over {\chi_{23}}}
~&=~{{4k+1}\over 2}{{q{\hat p}p^{\prime 2}}\over {p{\hat q}q^{\prime 2}}}
                                  {{\rho_{23}}\over {\chi_{23}}}
\cr
~&=~{{\rho_{12}}\over {\chi_{12}}}\cdot {{((4k+1)\rho_{23}{{p^{\prime}}
\over 2})^2}\over {(\chi_{23}q^{\prime 2})^2}}
\cr}
\eqno(36)
$$
Separating the square factors in ${{a_2}\over {b_2}}={{{\hat a}_2}\over
{{\hat b}_2}}\cdot {{{\hat p}^{\prime 2}}\over {{\hat q}^{\prime 2}}}$, 
it follows that
$$\eqalign{{{\hat a}_2}\over {{\hat b}_2}}{{((4k+1){p\over 2}\hat q)}\over
{(q{\hat p})}}~=~{{((4k+1){{p^\prime}\over 2} {\hat q}^\prime)^2}\over 
{(q^\prime {\hat p}^\prime)^2}}
\eqno(37)
$$
Either there is an overlap between the prime factors of $(4k+1){p\over 2}$ 
and ${\hat q}$ or ${\hat a}_2=(4k+1){p\over 2} 
{\hat q}=(4k+1){{p^\prime}\over 2} {\hat q}^\prime$, and similarly, 
there is either an overlap between the prime factors of $q$ and 
${\hat p}$ or ${\hat b}_2=q {\hat p}=q^\prime 
{\hat p}^\prime$.  Removing any overlap, then the remaining square factors 
can be separated in $a_2$ and $b_2$ obtaining the form 
${{{\hat a}_2}\over {{\hat b}_2}}$ for the square-free part of the 
ratio ${{a_2}\over {b_2}}$.  The equalities containing ${\hat a}_2$ and 
${\hat b}_2$ imply that ${\hat p} > {\hat p}^\prime \ge p^\prime > p$ and 
${\hat q} > {\hat q}^\prime \ge q^\prime > q$.  By interchanging the roles 
of $a_2,~b_2$ and $a_3,~b_3$ in the above argument, the inequalities 
$p> {\hat p}$ and $q > {\hat q}$ can be derived, implying a contradiction.  
Thus, when $\ell=3$, it should not be possible to find coefficients 
$\{a_i\}$ and $\{b_i\}$ satisfying equation (25) such that 
${{b_{13}}\over {a_{13}}}{{a_2}\over {b_2}}$ is ${{4k+1}\over 2}$ times
the square of a rational number.  The validity of this result is confirmed
by the non-existence of odd perfect numbers with four different prime 
factors.  

A variation of the standard induction argument can be used to show that 
there cannot be different odd perfect numbers with prime decompositions 
$(4k+1)^{4m+1}\prod_{i=1}^{\ell-1}q_i^{2\alpha_i}$ and
\hfil\break 
$(4k^\prime+1)^{4m^\prime+1}\prod_{i=1}^\ell~q_i^{\prime 2\alpha^\prime_i}$.  

When $\ell$ is odd, 
$${{q_1^{2\alpha_1+1}-1}\over {q_1-1}}...{{q_{\ell-1}^{m_{\ell-1}}-1}\over
{q_{\ell-1}-1}}~=~{{b_1...b_{\ell-1}}\over {a_1...a_{\ell-1}}}
\left({{(4k+1)^{4m+2}-1}\over {4k}}\right)^{\ell-1}~=~
{{b_1...b_{\ell-1}}\over {a_1...a_{\ell-1}}}\cdot \square
\eqno(38)
$$
rationality of square root of the product of repunits with 
$\ell-1$ prime bases $\{q_i,i=1,...,\ell-1\}$ would require
$${{b_1...b_\ell}\over {a_1...a_\ell}}~=~2(4k+1)\rho_\ell {{b_\ell}
\over {a_\ell}}\cdot \square
\eqno(39)
$$
and
$${{b_1...b_\ell}\over {a_1...a_\ell}}~=~2(4k+1)\rho_\ell\cdot \square
\eqno(40)
$$
Since the values $q_\ell=3$ and $\alpha_\ell=2$ can be excluded from the
product of repunits,  $\rho_\ell$ is odd and does not 
equal 1, so that ${{b_1...b_\ell}\over {a_1...a_\ell}}\ne 2(4k+1)\square$.  
The square root of the product of repunits with $\ell$ prime bases 
$\{q_i,i=1,...,\ell\}$ is therefore not rational. 

When $\ell$ is even,
$${{q_1^{2\alpha_1+1}-1}\over {q_1-1}}...{{q_{\ell-1}^{2\alpha_{\ell-1}+1}-1}
\over {q_{\ell-1}-1}}~=~{{b_1...b_{\ell-1}}\over {a_1...a_{\ell-1}}}
\left({{(4k+1)^{4m+2}-1}\over {4k}}\right)\cdot \square
\eqno(41)
$$
so that rationality of the square root expression with $\ell-1$ primes
$\{q_i,~i=1,...,\ell-1\}$ requires 
$${{b_1...b_\ell}\over {a_1...a_\ell}}~=~2(4k+1)\rho_\ell~\cdot~
\left({{(4k+1)^{4m+2}-1}\over {4k}}\right)~\cdot~\square
\eqno(42)
$$
Again, since $\rho_\ell\ne 1$, equation (42) implies that 
${{b_1...b_\ell}\over {a_1...a_\ell}}\left({{(4k+1)^{4m+2}-1}\over {4k}}\right)
~\ne~2(4k+1)\square$ or equivalently that the square root expression
with $\ell$ primes $\{q_i,~i=1,...,\ell\}$ is not rational.

The proof can be continued for $\ell > 3$ by assuming that there 
do not exist any odd primes $q_1,...,~q_{\ell-1}$ and 
$4k+1$ such that ${\sqrt {2(4k+1)}}\left[{{q_1^{2\alpha_1+1}-1}
\over {q_1-1}}...{{q_{\ell}^{2\alpha_{\ell}+1}-1}\over {q_{\ell-1}-1}}
\right]^{1\over 2}\left({{(4k+1)^{4m+2}-1}\over {4k}}\right)^{1\over 2}$ is 
rational and proving that the same property is valid when $\ell$ odd 
primes $q_1,...,~q_\ell$ arise in the prime decomposition of the integer $N$.

If $\ell$ is odd, $\left({{(4k+1)^{4m+2}-1}\over {4k}}
\right)^{{(\ell+1)}\over 2}$ is integer, and non-existence of odd perfect
numbers of the form $(4k+1)^{4m+1}q_1^{2\alpha_1}...
q_{\ell-1}^{2\alpha_{\ell-1}}$ is equivalent to the condition
${{b_1..b_{\ell}}\over {a_1...a_{\ell}}}\ne 2(4k+1) \square$.  Since
$$\eqalign{2(4k+1){{q_1^{m_1}-1}\over {q_1-1}}...{{q_{l-1}^{m_{\ell-1}-1}-1}
\over {q_{\ell-1}-1}}\left({{(4k+1)^{4m+2}-1}\over {4k}}\right)
~=&~2(4k+1){{b_1...b_{\ell-1}}\over {a_1...a_{\ell-1}}}\cdot 
\cr
&~~~~\left({{(4k+1)^{4m+2}-1}\over {4k}}\right)\cdot \square
\cr}
\eqno(43)
$$
Since the irrationality of the square root expression is assumed to hold
generally for $\ell-1$ odd primes $\{q_i\}$ and any value of $4k+1$, the 
effect of the inclusion of another prime $q_\ell$ can be deduced.  
Thus, given an arbitrary set of $\ell$ odd primes, $q_1,~...,~q_\ell$ and 
some prime of the form $4k+1$, irrationality of the square root of 
expression (43) implies that
$${{b_1...b_{\ell-1}}\over {a_1...a_{\ell-1}}}~\ne~
2(4k+1)~\left({{(4k+1)^{4m+2}-1}\over {4k}}\right)\cdot \square
\eqno(44)
$$
However, by equation (25), 
$\left({{(4k+1)^{4m+2}-1}\over {4k}}\right)={{a_\ell}\over {b_\ell}}
{{q_\ell^{2\alpha_\ell+1}-1}\over {q_\ell-1}}$, and if
${{q_\ell^{2\alpha_\ell+1}-1}\over {q_\ell-1}}\equiv \rho_\ell \chi_\ell^2$,
separating the square-free factors from the factors with even exponents.
it follows that
$$\eqalign{ {{b_1...b_{\ell-1}}\over {a_1...a_{\ell-1}}}~&\ne
~2(4k+1)~\rho_\ell {{a_\ell}\over {b_\ell}}\cdot \square
\cr
{{b_1...b_\ell}\over {a_1...a_\ell}}~&\ne~2(4k+1)~\rho_\ell \cdot \square
\cr}
\eqno(45)
$$

The form of the relation (45) is valid for arbitrary values of 
${{b_\ell}\over {a_\ell}}$, but the choice of $\rho_\ell$ is specific to 
the repunit
${{{q_\ell}^{2\alpha_\ell+1}-1}\over {q_\ell-1}}$.  Since 
${{{q_\ell}^{2\alpha_\ell+1}-1}\over {q_\ell-1}}$ is the square of an 
integer only when $q_\ell=3,~\alpha_\ell=2$, it is preferable to 
represent the rationality condition for $\ell-1$ and $\ell$ primes 
$\{q_i\}$ as
$$\eqalign{{{b_1...b_{\ell-1}}\over {a_1...a_{\ell-1}}}~&=~2 (4k+1) 
~\omega_{\ell-1} \rho_\ell {{a_\ell}\over {b_\ell}}\cdot \square
\cr
{{b_1...b_\ell}\over {a_1...a_\ell}}~&=~2(4k+1)~\omega_\ell\cdot\square
\cr}
\eqno(46)
$$
when $\ell$ is odd.  Irrationality of the square root expression for
$\ell-1$ primes $\{q_i,~i=1,...\ell-1\}$, which requires that 
$\omega_{\ell-1}\ne 1$ is a square-free 
integer, implies irrationality for $\ell$ primes $\{q_i,i=1,...,\ell\}$
if $\omega_{\ell-1}\rho_\ell=\omega_\ell\ne 1$ is square-free.

When $\ell$ is even, odd perfect numbers of the form
$(4k+1)^{4m+1} q_1^{2\alpha_1}...q_{\ell-1}^{2\alpha_{\ell-1}}$ do not
exist if ${{b_1...b_{\ell-1}}\over {a_1...a_{\ell-1}}}\ne 2 (4k+1)\cdot
\square$.  Then ${{b_1...b_\ell}\over {a_1...a_\ell}}\cdot 
\left({{(4k+1)^{4m+2}-1}\over {4k}}\right)\ne 2 (4k+1) \rho_\ell \cdot 
\square$.  Irrationality of the square root expression with $\ell-1$ primes 
$\{q_i,i=1...,\ell-1\}$ also can be represented as
as ${{b_1...{b_{\ell-1}}}\over {a_1...a_{\ell-1}}} = 2(4k+1)\omega_{\ell-1} 
\cdot \square$ where $\omega_{\ell-1}\ne 1$ is a square-free number.
Consequently, ${{b_1...b_\ell}\over {a_1...a_\ell}}= 2 (4k+1)\omega_{\ell-1} 
{{b_\ell}\over {a_\ell}} \cdot \square$.  Since irrationality of the square 
root expression with $\ell$ primes $\{q_i, i=1,...\ell\}$ would equivalent to
$$\eqalign{{{b_1...b_\ell}\over {a_1...a_\ell}}~\left({{(4k+1)^{4m+2}-1}
\over {4k}}\right)~&=~2(4k+1)~\omega_\ell \cdot \square
\cr
{{b_1...b_\ell}\over {a_1...a_\ell}}~&=~2(4k+1)~\omega_\ell \rho_\ell
                                                  {{a_\ell}\over {b_\ell}}
                                                        \cdot \square
\cr}
\eqno(47)
$$
this again can be achieved if $\omega_{\ell-1}\rho_\ell=\omega_\ell\cdot
\square$.  

For any prime divisor $p$
$$\eqalign{v_p(\omega_{\ell-1})~&=~\sum_{i=1}^{\ell-1}~
\biggl[v_p\left({{q_i^{e_i}-1}\over {q_i-1}}\right)
\delta\left({{n_i}\over {e_i}}-\left[{{n_i}\over {e_i}}\right]\right)
~+~v_p(n_i)\biggr]
\cr
&~~~~~~~~~~~~~~~~~~+~v_p\left({{(4k+1)^{4m+2}-1}\over {4k}}\right)~~~~~~
~~~~~~(mod~2)
\cr
v_p(\omega_\ell)~&=~\sum_{i=1}^\ell~\biggl[v_p\left({{q_i^{e_i}-1}\over 
{q_i-1}}\right)\delta\left({{n_i}\over {e_i}}-\left[{{n_i}\over 
{e_i}}\right]\right)~+~v_p(n_i)\biggr]
\cr
&~~~~~~~~~~~~~~~~~~+~v_p\left({{(4k+1)^{4m+2}-1}\over {4k}}\right)
~~~~~~~~~~~~
(mod~2)
\cr}
\eqno(48)
$$
where $e_i=ord_p(q_i)$.
It follows that
$$v_p(\omega_\ell)~=~v_p(\omega_{\ell-1})~+~v_p\left({{q_\ell^{e_\ell}-1}
\over {q_\ell-1}}\right)~+~v_p(n_\ell)
\eqno(49)
$$

Suppose that $p$ is one of the extra prime divisors so that 
$v_p(\omega_{\ell-1})=1$.  If $e_\ell\nmid n_\ell$ or $p\nmid n_\ell$,
then $p\not\bigg\vert~{{{q_\ell}^{n_\ell}-1}\over {q_\ell-1}}$ and 
$v_p(\omega_\ell)=1$.  

If $p^h \bigg\vert\bigg\vert {{q_\ell^{e_\ell}-1}\over {q_\ell-1}}$, and 
$p$ is a primitive prime divisor of this repunit, then $v_p(n_\ell)=0$ and 
$v_p(\omega_\ell)=1+h~(mod~2)$.  Since $v_p(\omega_\ell)=0~(mod~2)$ if 
$h=1$, it would be the next category of prime divisors, with the property 
$v_p\left({{q_\ell^{e_\ell}-1}\over {q_\ell-1}}\right)=2$ or equivalently
$Q_{q_\ell}\equiv 0~(mod~p)$, which contributes non-trivially to a square-free 
coefficient $\omega_\ell$. 
  
Since it has been assumed that the square root expression with $\ell-1$ 
primes $\{q_i,i=1,...,\ell-1\}$ is irrational, there is either an 
unmatched primitive divisor or an imprimitive divisor in the product 
$\prod_{i=1}^{\ell-1}~{{q_i^{2\alpha_i+1}-1}\over {q_i-1}}
~\cdot~{{(4k+1)^{4m+2}-1}\over {4k}}$.   Suppose that the extra 
prime divisor ${\hat p}_j$ is a factor of the repunit ${{q_j^{2\alpha_j+1}-1}
\over {q_j-1}}$.  By equation (27), 
$${{q_j^{2\alpha_j+1}-1}\over {q_j-1}}~=~\rho_j~\chi_j^2
~=~{{b_{j\ell}}\over {a_{j\ell}}} {{q_\ell^{2\alpha_\ell+1}-1}\over 
{q_\ell-1}}~=~{{b_{j\ell}}\over {a_{j\ell}}}~\rho_\ell~\chi_\ell^2
\eqno(50)
$$
so that $\rho_j \rho_\ell={{b_{j\ell}}\over {a_{j\ell}}}\cdot \square$.

To proceed further, it is first useful to choose the exponent $2\alpha_\ell+1$
to be equal to $2\alpha_j+1$.  If $p\vert (q_j-1)$, $p^{{\hat h}_j}\vert 
(2\alpha_j+1)$, $p\vert (q_\ell-1)$, $p^{{\hat h}_\ell}\vert 
(2\alpha_\ell+1)$, then $p^{{\hat h}_j}\bigg\vert 
{{q_j^{2\alpha_j+1}-1}\over {q_j-1}}$ and $p^{{\hat h}_\ell}\bigg\vert 
{{q_\ell^{2\alpha_\ell+1}-1}\over {q_\ell-1}}$  
When $\alpha_j=\alpha_\ell$, $p^{h_j}=p^{{\hat h}_j}=p^{{\hat h}_\ell}
=p^{h_\ell}$, where $h_j$ and $h_\ell$ denote the exponents of $p$
exactly dividing the repunits with bases $q_j$ and $q_\ell$ respectively,
so that this prime divisor will be absorbed into the square factors.

If $p\vert (q_j-1)$ and $p\nmid (q_\ell-1)$, then $h_j={\hat h}_j$ and 
$h_\ell={\hat h}_\ell+v_p\left({{q_\ell^{e_\ell}-1}\over {q_\ell-1}}
\right)$.  Since ${\hat h}_j={\hat h}_\ell$ when $\alpha_j=\alpha_\ell$,
$h_\ell=h_j+v_p\left({{q_\ell^{e_\ell}-1}\over {q_\ell-1}}\right)$.
Matching of the prime factors in the two repunits would require
$v_p\left({{q_\ell^{e_\ell}-1}\over {q_\ell-1}}\right)=0~(mod~2)$.
Because $p\vert (q_\ell^{e_\ell}-1)$, the minimum value of this exponent
is 2, implying that $Q_{q_\ell}\equiv 0~(mod~p)$.  Conversely, if
$Q_{q_\ell}\not\equiv 0~(mod~p)$ or $Q_{q_\ell}\equiv 0~
(mod~p^{h_\ell^\prime-1})$, $Q_{q_\ell}\not\equiv 0~(mod~p^{h_\ell^\prime})$,
where $h_\ell^\prime$ is odd, the prime divisor $p$ in the product of the
two repunits cannot be entirely absorbed into the square factors.    
Similar conclusions hold when $p\nmid (q_j-1)$ and $p\vert (q_\ell-1)$.

Let $p$ be an imprimitive prime divisor such that $p\nmid (q_j-1)$ and 
$p\nmid (q_\ell-1)$, then $v_p\left({{q_j^{n_j}-1}\over {q_j-1}}
\right)=v_p(q_j^{n_j}-1)$ and $v_p\left({{q_\ell^{n_\ell}-1}
\over {q_\ell-1}}\right)=v_p(q_\ell^{n_\ell}-1)$.  If $p^h\vert n_\ell$, and 
$n_j=n_\ell$, then $h_j=h_\ell=h$, again implying that the prime divisor 
can be absorbed into the square factors.

The arithmetic primitive factors of ${{q_j^{n_j}-1}
\over {q_j-1}}$ and ${{q_\ell^{n_\ell}-1}\over {q_\ell-1}}$, 
${{\Phi_{n_j}(q_j)}\over {p_j}}$ and ${{\Phi_{n_\ell}(q_\ell)}\over {p_\ell}}$ 
respectively, are different when $n_j=n_\ell$, except possibly for 
solutions generated by the prime equation ${{q_\ell^n-1}\over {q_j^n-1}}=p$ 
required when either $p_j=gcd(n_j,\Phi_{n_j}(q_j))$ or 
$p_\ell=gcd(n_\ell,\Phi_{n_\ell}(q_\ell))$ equals 1.  The algebraic 
primitive factors $\Phi_{n_j}(q_j)$ and $\Phi_{n_\ell}(q_\ell)$ necessarily 
will be different if $n_j=n_\ell$.  Consider a prime divisor $p^\prime$ of 
the arithmetic primitive factors which is raised to a different power in 
${{\Phi_{n_j}(q_j)}\over {p_j}}$ and ${{\Phi_{n_\ell}(q_\ell)}\over 
{p_\ell}}$.  If this prime is the only factor with this property, then 
${{q_\ell^{n_\ell}-1}\over {q_j^{n_j}-1}}={{q_\ell^{n_\ell}-1}\over
{q_j^{n_\ell}-1}}=(p^\prime)^{h_\ell-h_j}$, and the non-existence
of solutions to this equation for $h_\ell-h_j\ge 2$ has been shown in 
$\S 4$.  

The error in the approximation is given by ${{q_\ell^{n_\ell}}\over 
{q_j^{n_\ell}}}\left[1-{1\over {q_\ell^{n_\ell}}}+{1\over {q_j^{n_j}}}
+{\cal O}\left({1\over {q_j^{n_\ell}q_\ell^{n_\ell}}}
\right)\right]$, and since $\bigg\vert {1\over {q_j^{n_\ell}}}-
{1\over {q_\ell^{n_\ell}}} \bigg\vert < min\left({1\over {q_j^{n_\ell}}},
{1\over {q_\ell^{n_\ell}}}\right)$, the error is less than 
${{q_\ell^{n_\ell}}\over {q_j^{n_\ell}(q_j^{n_\ell}-1)}}\simeq
{{q_\ell^{n_\ell}}\over {q_j^{2n_\ell}}}$.   Given a rational 
number ${a\over b}$, the inequality $\bigg\vert {a\over b}-{{z_2}\over {z_1}}
\bigg\vert < {1\over {z_1^2}}$ has a finite number of solutions satisfying 
$z_1 < b$, $gcd(z_1,z_2)=1$ [39].  In particular, solutions to 
$$\bigg\vert {{{q_\ell}^{n_\ell}}\over {q_j^{n_\ell}}}-{{z_2}\over {z_1}}
\bigg\vert~=~\bigg\vert {{q_\ell^{n_\ell}}\over {q_j^{n_\ell}}}-{{y_2^2}
\over {y_1^2}}\bigg\vert~<~{1\over {y_1^4}}
\eqno(51)
$$
will be constrained by the inequality $y_1 < q_j^{{n_\ell}\over 2}$.
The condition $\bigg\vert {{q_\ell^{n_\ell}}\over {q_j^{n_\ell}}}
-{{y_2^2}\over {y_1^2}}\bigg\vert < {{q_\ell^{n_\ell}}\over {q_j^{2n_\ell}}}$
satisfied when ${{q_j^{{n_\ell}\over 2}}\over {q_\ell^{{n_\ell}\over 4}}}< y_1 
< q_j^{{n_\ell}\over 2}$.  
     
Since it has been established that square classes of the repunits ${{q^n-1}
\over {q-1}}$ consist of only one element [40], it follows that 
$(q_\ell^{n_{\ell_1}}-1)(q_\ell^{n_{\ell_2}}-1)=(\kappa^\prime)^2 (q_\ell-1)^2
(y_1^\prime)^2 (y_2^\prime)^2$ and there is only one representative from each 
sequence $\{q_j^{n_j}-1,~n_j\in {\Bbb Z}\}$, $\{q_\ell^{n_\ell}-1,~n_\ell\in 
{\Bbb Z}\}$  which has a specified square-free factor $\kappa$.  Thus, 
${{q_\ell^{n_\ell}-1}\over {q_j^{n_\ell}-1}}\ne {{y_2^2}\over {y_1^2}}$ unless 
$n_{q_j}(\kappa)$ coincides with $n_{q_\ell}(\kappa)$.  If   
$q_\ell^{n_\ell}-1=\kappa (y_2^\prime)^2$ and $q_j^{n_\ell}-1=\kappa 
(y_1^\prime)^2$, and ${{y_2^2}\over {y_1^2}}$ is the irreducible form
of ${{(y_2^\prime)^2}\over {(y_1^\prime)^2}}$, it follows that
$y_1 < {{q_j^{{n_\ell}\over 2}}\over {\sqrt {\kappa {\hat \kappa}^2}}}$, 
where ${\hat \kappa}=gcd (y_1^\prime,y_2^\prime)$.  Both inequalities for
$y_1$ cannot be satisfied if $q_\ell^{{n_\ell}\over 4} < {\sqrt {\kappa 
{\hat \kappa}^2}}$ or equivalently $q_\ell^{{n_\ell}\over 2} < 
gcd (q_j^{n_\ell}-1, q_\ell^{n_\ell}-1)$.  When the pair of primes 
$(q_j,q_\ell)$ satisfies the last inequality, the prime divisors in 
$\rho_j$ and $\rho_\ell$ do not match and the product
of the repunits ${{q_j^{n_j}-1}\over {q_j-1}}$ and ${{q_\ell^{n_\ell}-1}
\over {q_\ell-1}}$, with $n_j\ne n_\ell$, is not a perfect square.

The number of solutions to the inequality $\vert ax^n-by^n\vert\le h$ when 
$x~\ge~\left({{2h}\over {a^{1-\rho}\alpha}}\right)^{1\over {{n\over 2}-1}}$
with $\alpha=\left({b\over a}\right)^{1\over n}$ does not exceed
$6+{1\over {ln {n\over 2}}}\left[29+ln\rho^{-1}+ln\left(1+{{ln~2h}
\over {ln~a}}\right)\right]$ [41].  Setting 
${{q_\ell^n-1}\over {q_j^n-1}}\simeq {{y_2^2}\over {y_1^2}}$, it 
follows that $y_1^2(q_\ell^n-1)\simeq
y_2^2(q_j^n-1)$ leading to consideration of the inequality
$\vert y_2^2q_j^n-y_1^2 q_\ell^n\vert \le \vert y_2^2- y_1^2\vert$.     
The constraint placed on $q_j$ is
$$q_j\ge \Biggl({{2\vert y_2^2-y_1^2\vert}\over {y_2^{2(1-\rho)}
\left({{y_2^2}\over {y_1^2}}\right)^{1\over n}}}\Biggr)^{1\over {{n\over 2}-1}}
\eqno(52)
$$
Since ${{q_j-1}\over {q_\ell-1}}\ge \left({{y_1^2}\over {y_2^2}}
\right)^{1\over n}\ge {{q_j}\over {q_\ell}}$, it is sufficient for $q_j$ to
satisfy the stronger constraint
$$q_j\ge \biggl(2{{q_\ell-1}\over {q_j-1}}
                    y_2^{2\rho}\biggr)^{1\over {{n\over 2}-1}}
\eqno(53)
$$
which is equivalent to an upper bound for $y_2$ of
$$y_2^2\ge q_j^{\rho^{-1}\left({n\over 2}-1\right)}\cdot
\left({1\over 2}{{q_j-1}\over {q_\ell-1}}\right)^{\rho^{-1}}
\eqno(54)
$$ 
This condition defines an allowable range of values for $y_2$ when $\rho\le
{1\over 2}$.  The number of solutions to the inequality is not greater than
$$\eqalign{6~+~{1\over {ln~({n\over 2})}}\biggl(29~+~ln~\rho^{-1}~+~ln 
&\left(1~+~{{ln~(2\vert y_2^2-y_1^2\vert)}\over {ln~y_2^2}}\right)\biggr)
\cr
&~~~~~~~~~~~~~\le~6~+~{{29~+~ln(2+ln2)~+~ln~\rho^{-1}}\over {ln~{n\over 2}} }
\cr}
\eqno(55)
$$

Any adjustment in $n_\ell$ will introduce additional prime divisors.  Either 
they shall be new prime factors of the exponent or primitive divisors 
[42]-[45].   If $n_\ell$ is multiplied by a prime factor ${\hat p}^{r_\ell}$, 
where $\hat p\vert \rho_\ell$, then the product $\rho_\ell {\hat p}^{r_\ell}$ 
will contain the power ${\hat p}_\ell^{1+r_\ell}$.  While the prime power can 
be absorbed into the product of square factors only when $r_\ell$ is odd,  
the repunit ${{q_\ell^{n_\ell{\hat p}^{r_\ell}}-1}\over {q_\ell-1}}$ has
extra primitive divisors, giving rise to a non-trivial $\omega_\ell$, 
implying irrationality of the square root expression with $\ell$ primes 
$\{q_i,~i=1,...\ell\}$.   Moreover, $gcd(\Phi_{{\hat p}^i}(q), 
\Phi_{{\hat p}^j}(q))=1$ when $i\ne j$ and $p\nmid (q-1)$, 
multiplication of the index by ${\hat p}^{r_\ell}$ will introduce new 
prime divisors through the decomposition of the repunit 
${{q_\ell^{n_\ell{\hat p}^{r_\ell}}-1}\over {q_\ell-1}}
=\prod_{{d\vert n_\ell{\hat p}^{r_\ell}}\atop {d>1}}~\Phi_d(q_\ell)$.

The abstract argument given for $\ell=3$ could also be extended to higher
values of $\ell$.  This approach would consist of the demonstration of
the property ${{b_1...b_\ell}\over {a_1...a_\ell}}\cdot \ne 2(4k+1)\cdot 
\square$ if $\ell$ is odd, and  ${{b_1...b_\ell}\over {a_1...a_\ell}}
\cdot \ne 2(4k+1)\left({{(4k+1)^{4m+2}-1}\over {4k}}\right)$ if $\ell$
is even, given that there are no sets of primes $\{q_i\}$ with less than 
$\ell$ elements satisfying the rationality condition.   It 
may be noted that since
$$\eqalign{{{b_1...b_\ell}\over {a_1...a_\ell}}~&=~\left({{b_{13}}\over 
{a_{13}}}{{a_2}\over {b_2}}\right)\left({{b_{46}}\over {a_{46}}}
{{a_5}\over {b_5}}\right)...\left({{b_{\ell-2,\ell}}\over {a_{\ell-2,\ell}}}
{{a_{\ell-1}}\over {b_{\ell-1}}}\right)\cdot 
\left({{b_2b_3b_5b_6...b_{\ell-1}b_\ell}\over {a_2a_3a_5a_6...
a_{\ell-1}a_\ell}}\right)^2
\cr
&~~~~~~~~~~~~~~~~~~~~~~~~~~~~~~~~~~~~~~~~~~~~~~~~~~~~~~~~~~~~~~~~~~~~~~
when~\ell\equiv~0~(mod~3)
\cr
{{b_1...b_\ell}\over {a_1...a_\ell}}~&=
~\left({{b_{13}}\over {a_{13}}}
{{a_2}\over {b_2}}\right)\left({{b_{46}}\over {a_{46}}}{{a_5}\over {b_5}}
\right)...\left({{b_{\ell-3,\ell-1}}\over {a_{\ell-3,\ell-1}}}
{{a_{\ell-2}}\over {b_{\ell-2}}}\right){{b_\ell}\over {a_\ell}}
\cdot \left({{b_2b_3b_5b_6...b_{\ell-2}b_{\ell-1}}\over {a_2a_3a_5a_6...
a_{\ell-2}a_{\ell-1}} }\right)^2
\cr
&~~~~~~~~~~~~~~~~~~~~~~~~~~~~~~~~~~~~~~~~~~~~~~~~~~~~~~~~~~~~~~~~~~~~~~when~
\ell\equiv~1~(mod~3)
\cr
{{b_1...b_\ell}\over {a_1...a_\ell}}~&=~\left({{b_{13}}\over {a_{13}}}
{{a_2}\over {b_2}}\right)\left({{b_{46}}\over {a_{46}}}{{a_5}\over {b_5}}
\right)...\left({{b_{\ell-4,\ell-2}}\over {a_{\ell-4,\ell-2}}}
{{a_{\ell-3}}\over {b_{\ell-3}}}\right)
{{b_{\ell-1}b_\ell}\over {a_{\ell-1}a_\ell}}
\cdot \left({{b_2b_3b_5b_6...b_{\ell-3}b_{\ell-2}}\over 
{a_2a_3a_5a_6...a_{\ell-3}a_{\ell-2}}}\right)^2
\cr
&~~~~~~~~~~~~~~~~~~~~~~~~~~~~~~~~~~~~~~~~~~~~~~~~~~~~~~~~~~~~~~~~~~~~~~
when~\ell\equiv~2~(mod~3)
\cr}
\eqno(56)
$$
\noindent
and ${{b_{13}}\over {a_{13}}}{{a_2}\over {b_2}}=2(4k+1){{{\bar \rho}_1}
\over {{\bar \chi}_1}}
\cdot \square,~...,~{{b_{\ell-k^\prime-2,\ell-k^\prime}}\over
{a_{\ell-k^\prime-2,\ell-k^\prime}}}=2(4k+1)
{{{\bar \rho}_{\left[{\ell \over 3}\right]}}\over 
{{\bar \chi}_{\left[{\ell \over 3}\right]}}} \cdot 
\square$, where $\ell\equiv k^\prime~(mod~3)$,
\hfil\break
$k^\prime=0,1,2$,~${\bar \rho}_1$,...,~${\bar \rho}_{\left[{\ell \over 3}
\right]}$,~${\bar \chi}_1$,...,~${\bar \chi}_{\left[{\ell \over 3}\right]}$
are square-free factors, the quotient will equal
\hfil\break 
$\left(2(4k+1)\right)^{\left[{\ell \over 3}\right]}
f_{k^\prime}{{{\bar \rho}_1}\over {{\bar \chi}_1}}...
{{{\bar \rho}_{\left[{\ell \over 3}\right]}}\over 
{{\bar \rho}_{\left[{\ell\over 3}\right]}}}\cdot\square$ 
with $f_0=1$, $f_1={{b_\ell}\over {a_\ell}}$ and 
$f_2={{b_{\ell-1}b_\ell}\over {a_{\ell-1}a_\ell}}$.  It has been
established that ${{b_\ell}\over {a_\ell}}\ne 2(4k+1)\cdot \square$
because there is no odd integer of the form $(4k+1)^{4m+1}
q_\ell^{2\alpha_\ell}$ which satisfies the relation ${{\sigma(N)}\over N}=2$.
${{b_{\ell-1}b_\ell}\over {a_{\ell-1}a_\ell}}\ne 2(4k+1)\left({{4k+1)^{4m+2}-1}
\over {4k}}\right)\cdot \square$ because of the non-existence of odd perfect
numbers of the type $(4k+1)^{4m+1}q_{\ell-1}^{2\alpha_{\ell-1}}
q_\ell^{2\alpha_\ell}$.  Setting ${{b_\ell}\over {a_\ell}}=2(4k+1)
{{{\hat \rho}_{\ell 1}}\over {{\hat \chi}_1}}\cdot \square$ and 
${{b_{\ell-1}b_\ell}\over {a_{\ell-1}a_\ell}}=2(4k+1){{{\hat \rho}_{\ell 2}}
\over {{\hat \chi}_{\ell 2}}}\left({{{(4k+1)}^{4m+2}-1}
\over {4k}}\right)\cdot \square$, it follows
that
$$\eqalign{{{b_1...b_\ell}\over {a_1...a_\ell}}~&=~(2(4k+1))^{\ell\over 3}
{{{\bar \rho}_1}\over {{\bar \chi}_1}}...{{{\bar \rho}_{\ell\over 3}}
\over {{\bar \chi}_{\ell \over 3}}}\cdot \square
~~~~~~~~~~~~~~~~~~~~~~~~~~~~~~~~~~~~~~~~~~~~~~\ell\equiv~0~(mod~3)
\cr
{{b_1...b_\ell}\over {a_1...a_\ell}}~&=~
(2(4k+1))^{\left[{\ell\over 3}\right]+1}
{{{\bar \rho}_1}\over {{\bar \chi}_1}}...{{{\bar \rho}_{\left[{\ell\over 3}
\right]}}\over {{\bar \chi}_{\left[{\ell \over 3}\right]}}}
\cdot {{{\hat \rho}_{\ell 1}}\over {{\hat \chi}_{\ell 1}}}
\cdot \square~~~~~~~~~~~~~~~~~~~~~~~~~~~~~~~~\ell\equiv~1~(mod~3)
\cr
{{b_1...b_\ell}\over {a_1...a_\ell}}&=
(2(4k+1))^{\left[{\ell\over 3}\right]+1}
{{{\bar \rho}_1}\over {{\bar \chi}_1}}...{{{\bar \rho}_{\left[{\ell \over 3}
\right]}}
\over {{\bar \chi}_{\left[{\ell \over 3}\right]}}}\cdot {{{\hat \rho}_{\ell 2}}
\over {{\hat \chi}_{\ell 2}}}
\left({{(4k+1)^{4m+2}-1}\over {4k}}\right)\cdot \square
~~~~~\ell\equiv~2~(mod~3)
\cr}
\eqno(57)
$$
and the coefficients $\{a_i,b_i\}$ will not satisfy the rationality condition
when the square-free factors ${\bar \rho}_1,...,
{\bar \rho}_{\left[{\ell \over 3}\right]}, {\hat \rho}_{\ell 1}, 
{\hat \rho}_{\ell 2}, {\bar \chi}_1,...,
{\bar \chi}_{\left[{\ell \over 3}\right]}, {\hat \chi}_{\ell 1},
{\hat \chi}_{\ell 1}$ have prime divisors other than $2$ and $4k+1$
which do not match to produce the square of a rational number.  

When $\ell$ is odd and greater than $5$, there always exists an odd integer 
$\ell_o$ and an even integer $\ell_e$ such that $\ell=3\ell_o+2\ell_e$, 
implying the following identity 
$${{b_1...b_\ell}\over {a_1...a_\ell}}=\left({{b_{13}}\over {a_{13}}}{{a_2}
\over {b_2}}\right)\left({{b_{46}}\over {a_{46}}}{{a_5}\over {b_5}}\right)
...\left({{b_{3\ell_o-2,3\ell_o}}\over {a_{3\ell_0-2,3\ell_o}}}{{a_{3\ell_o-1}}
\over {b_{3\ell_o-1}}}\right)\left({{b_{3\ell_0+1}b_{3\ell_0+2}}\over
{a_{3\ell_0+1}a_{3\ell_0+2}}}\right)...\left({{b_{\ell-1}b_\ell}\over
{a_{\ell-1}a_\ell}}\right)\cdot \square
\eqno(58)
$$
Consequently,
$$\eqalign{{{b_1...b_\ell}\over {a_1...a_\ell}}~&=~(2(4k+1))^{\ell_0+\ell_e}
~\cdot~{{{\bar \rho}_1}\over {{\bar \chi}_1}}...{{{\bar \rho}_{\ell_o}}\over
{{\bar \chi}_{\ell_o}}}...{{{\hat \rho}_{\ell-2\ell_e+2,2}}\over
{{\hat \chi}_{\ell-2\ell_e+1,2}}}...{{{\hat \rho}_{\ell 2}}\over 
{{\hat \chi}_{\ell 2}}}~\cdot~\left({{(4k+1)^{4m+2}-1}\over {4k}}
\right)^{\ell_e}~\cdot~\square
\cr
~&=~2(4k+1)~\cdot~{{{\bar \rho}_1}\over {{\bar \chi}_1}}...
{{{\bar \rho}_{\ell_o}}\over {{\bar \chi}_{\ell_o}}}...
{{{\hat \rho}_{\ell-2\ell_e+2,2}}\over {{\hat \chi}_{\ell-2\ell_e+2,2}}}...
{{{\hat \rho}_{\ell 2}}\over {{\hat \chi}_{\ell 2}}}~\cdot~\square
\cr}
\eqno(59)
$$
Regardless of the factors of $2$ and $4k+1$, the coefficients $\{a_i,b_i\}$ 
will produce an irrational square root expression (26) for odd $\ell$
if the product of fractions $\prod_{i=1}^{\ell_o}
~{{{\bar \rho}_i}\over {{\bar \chi}_i}}~\prod_{j=1}^{\ell_e}
~{{{\hat \rho}_{\ell-2j+2,2}}\over {{\hat \chi}_{\ell-2j+2,2}}}$
is not the square of a rational number.

If $\ell$ is even and greater than $4$, there always exists an odd integer 
$\ell_o$ and an even integer $\ell_e$ such that $\ell=2\ell_o+3\ell_e$.  
From the identity
$${{b_1...b_\ell}\over {a_1...a_\ell}}~=~\left({{b_1b_2}\over {a_1a_2}}\right)
...\left({{b_{2\ell_o-1}b_{2\ell_o}}\over {a_{2\ell_o-1}a_{2\ell_o}}}\right)
...\left({{b_{2\ell_o+1,2\ell_0+3}}\over {a_{2\ell_o+1}a_{2\ell_o+3}}}
{{a_{2\ell_o+2}}\over {b_{2\ell_o+2}}}\right)...\left({{b_{\ell-2,\ell}}\over
{a_{\ell-2,\ell}}}{{a_{\ell-1}}\over {b_{\ell-1}}}\right)\cdot \square
\eqno(60)
$$
it follows that   
$$\eqalign{{{b_1...b_\ell}\over {a_1...a_\ell}}~&=~(2(4k+1))^{\ell_o+\ell_e}
\left({{(4k+1)^{4m+2}-1}\over {4k}}\right)^{\ell_o}~\cdot~
{{{\hat \rho}_{22}}\over {{\hat \chi}_{22}}}...{{{\hat \rho}_{2\ell_o,2}}
\over {{\hat \chi}_{2\ell_o,2}}} {{{\bar \rho}_{\ell-3\ell_e+1}}\over
{{\bar \chi}_{\ell-3\ell_e+1}}}...{{{\bar \rho}_{\ell-2}}\over
{{\bar \chi}_{\ell-2}}}~\cdot~\square
\cr
~&=~2(4k+1)~\cdot~{{{\hat \rho}_{22}}\over {{\hat \chi}_{22}}}...
{{{\hat \rho}_{2\ell_o,2}}\over {{\hat \chi}_{2\ell_o,2}}}
{{{\bar \rho}_{\ell-3\ell_e+1}}\over {{\bar \chi}_{\ell-3\ell_e+1}}}
...{{{\bar \rho}_{\ell-2}}\over {{\bar \chi}_{\ell-2}}}
~\cdot~\left({{(4k+1)^{4m+2}-1}\over {4k}}\right)~\cdot~\square
\cr}
\eqno(61)
$$
Again, the factors of $2$ and $4k+1$ are not relevant, and the coefficients
$\{a_i,b_i\}$ give rise to an irrational square root expression (26)
for even $\ell$ if $\prod_{i=1}^{\ell_o}~{{{\hat \rho}_{2i,2}}\over
{{\hat \chi}_{2i,2}}}~\prod_{j=1}^{\ell_e}~{{{\bar \rho}_{\ell-3j+1}}
\over {{\bar \chi}_{\ell-3j+1}}}$ is not the square of a rational number.
\qed
\vskip 10pt
\centerline{\bf 7. Conclusion}

The rationality condition provides an analytic method for investigating the
existence of odd perfect numbers.  The aim of this approach then becomes
the proof of the existence of an unmatched prime divisor in the product
of the repunits, since the square root of any such divisor would be 
irrational, contrary to the condition for the existence of an odd perfect 
number.  An upper bound for the density of odd integers greater than 
$10^{300}$, in an interval of fixed length, which could satisfy 
${{\sigma(N)}\over N}=2$, may be found by considering the square root 
expression containing the product of repunits, combining the estimate of the 
density of square-full numbers in this range with the probability of an 
integer being expressible as the product of repunits with prime bases 
multiplied by $2(4k+1)$.  Repunits form a special class of Lucas sequences, 
and the properties of primitive and imprimitive prime divisors of these 
sequences can be used to determine the powers of primes dividing the product 
of repunits.  A comparison of the divisibility properties of Lucas sequences 
$U_n(q+1,q)$ with different values of $q$ has been undertaken in $\S 4$.   
Specifically, the arithmetic primitive factors of these repunits,
products of the primitive prime power divisors, can be compared for different
values of the prime basis, and it has been shown that they could only
be equal if the indices of Lucas sequences differ, except possibly for pairs 
of divisors $\left(\Phi_n(q_i), {{\Phi_n(q_j)}\over {p_j}}\right)$
generated by the prime equation ${{q_j^n-1}\over {q_i^n-1}}=p$. In the
second theorem, non-existence of the odd perfect numbers for a large set of 
primes $\{q_i,i=1,...,\ell;~4k+1\}$, exponents $\{2\alpha_i+1,i=1,...,\ell;~
4m+1\}$ and values of $\ell$ using the method of induction adapted to 
the coefficients $\{a_i,b_i\}$ in the product of repunits.  An abstract
argument is given for the non-existence of coefficients satisfying the
rationality condition when $\ell=3$ and then various results are proven
for $\ell>3$ by using the properties of prime divisors of product of two
repunits, ${{q_j^{2\alpha_j+1}-1}\over {q_j-1}}$ and 
${{q_\ell^{2\alpha_\ell+1}-1}\over {q_\ell-1}}$,
belonging to each of the four categories: (i) $p\vert (q_j-1)$, 
$p\vert (q_\ell-1)$ (ii) $p\vert (q_j-1)$, $p\nmid (q_\ell-1)$
(iii) $p\nmid (q_j-1)$, $p\vert (q_\ell-1)$ (iv) $p\nmid (q_j-1)$,
$p\nmid (q_\ell-1)$.   Irrationality of the square root expression for any 
set of $\ell-1$ primes $\{q_i,~i=1,...,\ell-1\}$ implies that each unmatched 
prime divisor in the product of repunits with bases 
$\{q_i,~i=1,...,\ell-1,4k+1\}$ can be associated with a single 
repunit, because factors of other repunits divisible by this prime 
contain powers of the prime with the exponent summing up to an even integer.
Supposing, for example, that the repunit containing this extra
prime divisor is ${{q_j^{2\alpha_j+1}-1}\over {q_j-1}}$.  The problem of 
determining whether this prime divisor remains unmatched, when an
additional prime $q_\ell$ in the decomposition of the odd integer $N$ is 
included, depends on the feasibility of matching the prime divisors of
each pair of repunits $(U_{n_j}(q_j+1,q_j),U_{n_\ell}(q_\ell+1,q_\ell))$
as $j$ takes all values in the range $\{1,2,...,\ell-1\}$ such that the 
repunit $U_{n_j}(q_j+1,q_j)$ contains an extra prime divisor.   
 
\vskip 20pt
\centerline{\bf Acknowledgements}
\noindent
I would like to acknowledge useful discussions with Prof. T. Gagen, who
suggested a reduction of equation (12) that is useful in the proof of the 
non-existence of odd perfect numbers with a specific condition imposed on 
the prime divisors. 
\vskip 20pt
\noindent
\centerline{\bf References}
\baselineskip=10pt
\parskip=7pt
\item{[1]} J. J. Sylvester, `On the divisors of the sum of a geometrical
series whose first term is unity and common ratio any positive or negative
integer', Nature ${\underline{87}}$ (1888) 417-418
\item{[2]} P. Hagis, Math. Comp. ${\underline{35}}$ (1980) 1027-1032  
\item{[3]} M. Kishore, Math. Comp. ${\underline{32}}$ (1978) 303-309
\item{[4]} L. Euler, Tractatus de Numerorum Doctrina, $\S 109$ in
${\underline{Opera~Omnia}}$ I, 5 (Genevae:
\hfil\break
Auctoritate et Impensis Societatis Scientarum Naturalium Helveticae, MCMXLIV)
\item{[5]} J. J. Sylvester, Comptes Rendus CVI (1888), pp. 403-405
\item{[6]} J. A. Ewell, Journal of Number Theory ${\underline{12}}$ 
(1980) 339-342
\item{[7]} T. Nagell, Norsk. Mat. Tidsskr. ${\underline{2}}$ (1920)
75-78
\hfil\break
T. Nagell, Mat. Fornings Skr. ${\underline{1}}$(3)(1921) 
\item{[8]} W. Ljunggren, Norsk Mat. Tidsskr. ${\underline{25}}$ (1943) 17-20
\item{[9]} P. Ribenboim, ${\underline{Catalan's~Conjecture: 
~Are~8~and~9~the~Only~Consecutive~Powers?'}}$ (Sydney: Academic Press, 1994)
\item{[10]} R. P. Brent, G. L. Cohen and H.J.J te Riele,
Math. Comp. ${\underline{57}}$ (1991) 857-868
\item{[11]} P. Hagis and G. L. Cohen, Math. Comp. ${\underline{67}}$ (1998)
1323-1330
\item{[12]} D. Iannucci, Math. Comp. ${\underline{68}}$ (1999) 1748 - 1760;
D. Iannucci, Math. Comp. ${\underline{69}}$ (2000) 867 - 879
\item{[13]} P. Hagis, Math. Comp. ${\underline{40}}$ (1983) 399-404
\item{[14]} M. Kishore, Mat. Comp. ${\underline{40}}$ (1983) 405-411
\item{[15]} N. Robbins, J. Reine Angew. Math. ${\underline{279}}$ (1975) 
14-21
\item{[16]} L. Somers, Fib. Quart. ${\underline{18}}$(4) (1980) 316-334
\item{[17]} D. H. Lehmer Ann. Math. ${\underline{31}}$ (1930) 419-448
\item{[18]} M. Hall, Bull. American Math. Soc. ${\underline{43}}$ (1937) 
78-80
\item{[19]} J. Brillhart, J. Toscania and P. Weinberger, 
${\underline{Computers~in~Number~Theory}}$, (London: Academic Press, 1971)
\item{[20]} S. Yates, ${\underline{Repunits~and~Repetends}}$ (Boynton Beach,
Florida: Star Publishing Co.Inc., 1982)
\item{[21]} J. J. Sylvester, Amer. J. Math. ${\underline{2}}$ (1879) 365
\item{[22]} G. D. Birkhoff and H. S. Vandiver, Annals of Math.,
${\underline{5}}$ (1904) 173-180
\item{[23]} L. E. Dickson, Amer. Math. Monthly ${\underline{16}}$ (1905)
86-89
\item{[24]} A. Schinzel, Acta Arithmetica ${\underline{15}}$ (1968) 49-70
\item{[25]} C. L. Stewart, Acta Arithmetica ${\underline{26}}$ (1975)
427-433
\item{[26]} C. L. Stewart, Proc. London Math. Soc. ${\underline{35}}$
(1977) 425-447
\item{[27]} C. L. Stewart, `Primitive Divisors of Lucas and Lehmer Numbers',
\hfil\break
${\underline{Transcendence~Theory:~Advances~and~Applications}}$ 
(London: Academic Press, 
\hfil\break
1977), pp. 79-92
\item{[28]} K. Motose, Math. J. Okayama Univ. ${\underline{35}}$ (1993) 
35-40; Math. J. Okayama Univ. ${\underline{37}}$ (1995) 27-36
\item{[29]} B. Richter, J. Reine Angew. Math. ${\underline{267}}$ (1974) 
77-89
\item{[30]} R. Ernvall and T. Mets{\"a}nkyl{\"a}, Math. Comp.   
${\underline{66}}$ (1997) 1353-1365
\item{[31]} J. W. L. Glaisher, Quarterly Journal of Pure and Applied 
Mathematics ${\underline{32}}$ (1901) 1-27, 240-251 
\item{[32]} L. E. Dickson, ${\underline{History~of~the~theory~of~numbers}}$ 
Vol. 1 (New York: Chelsea, 1966)
\item{[33]} K. Hensel `Theorie der algebraischen Zahlen' (Leipzig: Teubner,
1908) 
\item{[34]} N. Koblitz, ${\underline{p-adic~numbers,~p-adic~analysis
~and~zeta-functions}}$ (New York: 
\hfil\break
Springer-Verlag, 1977), pp.16-18
\item{[35]} W. Johnson, J. Reine Angew. Math. ${\underline{292}}$ (1977) 
196-200
\item{[36]} R. D. Fray, Duke Math. J. ${\underline{34}}$ (1967) 467-480
\item{[37]} R. Balasubramanian and T. N. Shorey, Math. Scand. 
${\underline{46}}$ (1980) 177-182
\item{[38]} T. N. Shorey, Hardy-Ramanujan Journal, Vol. 7 (1984) 1-10
\item{[39]} T. Nagell, ${\underline{Introduction~to~Number~Theory}}$
(New York: Chelsea Publishing Company, 1964)
\item{[40]} P. Ribenboim, J. Sichuan Univ., Vol. 26, Special Issue (1989)
196-199 
\item{[41]} J. Mueller, Quart. J. Math. Oxford ${\underline{38}}$ (1987)
503-513
\item{[42]} A. S. Bang, Taltheoretiske Undersogelser, Tidskrifft for 
Math. ${\underline{5}}$ (1886) 70-80; 130-137
\item{[43]} K. Zsigmondy, Zur Theorie der Potenzreste, Monatsh. fur Math. 
${\underline{3}}$ (1892) 265-284
\item{[44]} G. D. Birkhoff and H. S. Vandiver, Annals of Math. 
${\underline{5}}$ (1904) 173-180
\item{[45]} R. D. Carmichael, Ann. of Math. ${\underline{15}}$ (1913) 30-70

\end